\documentclass[11pt]{article}

\usepackage{amsfonts}
\usepackage{mathrsfs} 
\usepackage{amsfonts,amssymb}
\usepackage{exscale}

\usepackage{comment}

\usepackage{amsgen,amsmath,amstext,amsbsy,amsopn,amsfonts,amssymb,pifont, epsfig}
 \usepackage{graphicx}

\newtheorem{thm}{Theorem}[section]

\newtheorem{prop}[thm]{Proposition}
\newtheorem{defn}[thm]{Definition}

\newtheorem{example}[thm]{Example}
\newtheorem{remark}[thm]{Remark}

\newcommand{\be}{\begin{equation}}
\newcommand{\ee}{\end{equation}}
\newcommand{\bea}{\begin{eqnarray}}
\newcommand{\eea}{\end{eqnarray}}
\newcommand{\beaa}{\begin{eqnarray*}}
\newcommand{\eeaa}{\end{eqnarray*}}
\newcommand{\bei}{\begin{itemize}}
\newcommand{\eei}{\end{itemize}}
\newcommand{\bee}{\begin{enumerate}}
\newcommand{\eee}{\end{enumerate}}
\newcommand{\bi}{\begin{itemize}}
\newcommand{\ei}{\end{itemize}}
\newcommand{\beq}{\begin{eqnarray*}}
\newcommand{\eeq}{\end{eqnarray*}}
\newcommand{\beqn}{\begin{eqnarray}}
\newcommand{\eeqn}{\end{eqnarray}}

\newcommand{\pf}{\noindent {\bf Proof:\ }}

\newcommand{\ignore}[1]{}{}

\def\EE{\mathcal{E}}
\def\FF{\mathcal{F}}
\def\LL{\mathcal L}
\def\DD{\mathcal D}

\def\H{{\mathbf H}}
\def\qed{{\hfill $\Box$ \bigskip}}

\def\eps{\varepsilon}

\def\<{\langle}
\def\>{\rangle}
\def\wh{\widehat}
\def\wt{\widetilde}

\def\P{{\mathbb P}}

\def\E{{\mathbb E}\,}

\def\R{{\mathbb R}}

\newcommand{\cC}{{\cal C}}

\numberwithin{equation}{section}

\begin{document}

\title{Markov processes with darning and their approximations}

\author{{Zhen-Qing Chen} \quad \hbox{and} \quad { Jun Peng}}

\date{February 6, 2017}
 
\maketitle

\begin{abstract}

In this paper, we  study darning of general symmetric Markov processes
  by shorting some parts of the state space into singletons.
A natural way to construct such processes is via Dirichlet forms restricted to
the function space whose members take constant values on these collapsing parts.
They include as a special case Brownian motion with darning, which has been studied
in details in \cite{Chen, CF, CFR}.
When the initial processes have discontinuous sample paths,
the processes constructed in this paper are the genuine  extensions of those studied in Chen and Fukushima \cite{CF}.
We further show that, up to a time change, these Markov processes with darning can be approximated
in the finite dimensional sense  
by introducing additional large intensity jumps among
 these compact sets to be collapsed into singletons to the original Markov processes. 
For diffusion processes, it is also possible to get, up to a time change, diffusions with darning by
increasing the conductance on these compact sets to infinity.
To accomplish these, we extend the semigroup characterization of Mosco convergence to closed symmetric  forms
whose domain of definition may not be dense in the $L^2$-space. The latter is of independent interest
and potentially useful to study convergence of Markov processes having different state spaces.
Indeed,  we show in Section \ref{S:5} of this paper that Brownian motion in a plane with a very thin flag pole 
can be approximated
by Brownian motion in the plane with a vertical cylinder whose horizontal motion on the cylinder 
is a circular Brownian motion   moving at  fast speed.  
\end{abstract}



\bigskip

\section{ Introduction }

K. Ito \cite{KI} introduced the notion of  Poisson point process of excursions around one point $a$ in the state space of a standard Markov process $X$. He was motivated by a giving systematic constructions of Markovian extensions of the absorbing diffusion process $X^0$ on the half line $(0,\infty)$ subject to Feller's general boundary conditions \cite{IM}. Ito had constructed the most general jump-in process from the exit boundary $0$ by using  Poisson point process of excursions.
 Recent study \cite{FT, CFY, CF} reveals that Ito's program works equally well in the study of   Markov processes transformed 
 by collapsing  certain compact subsets of the state space  into singletons. These processes are called Markov processes with darning in \cite{CF}.
 (When the underlying process is a Markov chain on a discrete state space, such a procedure of collapsing subsets of state space is also called
 shorting in some literature.) 
However, in order to use excursion theory, it is assumed in \cite{FT, CFY, CF} that the original Markov 
process enters these compact subsets in a continuous way. This condition is automatically satisfied
for diffusion processes but not for general symmetric Markov processes that may have discontinuous trajectories. 

The purpose of this paper is two-folds. First, we extend the notion and construction of
Markov processes with darning to any symmetric Markov process, without assuming that 
the processes enter the compact subsets to be collapsed  in a continuous way.
In this generality, we can no longer use Poisson point process of excursions for the construction.
We will use instead  a Dirichlet form approach, which turns out to be quite effective. 
The second goal is to investigate  approximation schemes 
for  general Markov processes with darning by more concrete processes, which can be used
for simulation. For this, we develop Mosco convergence of closed symmetric    forms whose domain
may not be dense in the underlying $L^2$-space. This is because due to the collapsing of
the compact holes, the domain of the Dirichlet form for the Markov process with darning
is not dense in the $L^2$-space on the original state space. We now describe the content of
this paper in some details. For basic definitions and properties of symmetric Dirichlet forms, 
we refer the reader to \cite{CF, FOT}.

Let $E$ be a locally compact separable metric space and $m$ a Radon measure on $E$ with full support.
Suppose $(\EE, \FF)$ is a regular  Dirichlet form on $L^2(E; m)$
 in the sense that $C_c(E)\cap \mathcal{F}$ is dense both in $C_c(E)$ with respect to the uniform norm in $\mathcal{F}$
 and with respect to the Hilbert norm $\sqrt{\mathcal{E}_1(u,u)}:=\sqrt{\mathcal{E}(u,u)+(u,u)_{L^2(E; m)}}$.
 Here and in the sequel, we use $:=$ as a way of definition and 
$C_c(E)$ is the space of continuous functions on $E$ with compact support.
Every $f$ in $\FF$ admits an $\EE$-quasi-continuous $m$-version, which is unique up to
an $\EE$-polar set. We always take such a quasi-continuous version for functions in $\FF$. 
There is an $m$-symmetric Hunt process $X$ on $E$ associated with $(\EE, \FF)$, which is unique up to an $\EE$-polar set. It is known that for any regular Dirichlet form $(\EE, \FF)$ on $L^2(E; m)$
it admits the following unique Beurling-Deny decomposition (see \cite{CF, FOT}):
$$
\EE (u, u)= \EE^c (u, u)+\frac12 \int_{E\times E} (u(x)-u(y))^2 J(dx, dy) + \int_E u(x)^2 \kappa (dx),
\quad u\in \FF, 
$$
where $\EE^c$ is a symmetric non-negative definite
 bilinear form on $\FF$ that satisfies strong local property,
where $J(dx, dy)$ is a $\sigma$-finite measure on $E\times E\setminus {\rm diagonal}$,
and $\kappa$ is a $\sigma$-finite smooth measure on $E$. 
The measures $J(dx, dy)$ and $\kappa (dx)$ are called the jumping measure and killing of
the process $X$ (or equivalently, of the Dirichlet form $(\EE, \FF)$). 
Indeed, if we use $(N(x, dy), H_t)$ to denote the L\'evy system of $X$,
where $N(x, dy)$ is a kernel on $E_\partial:=E\cup\{\partial\}$ and $t\mapsto H_t$ 
is a positive continuous additive functional (PCAF) of $X$, then 
$$
J(dx, dy)=N(x, dy) \mu_H (dx) \qquad \hbox{and}  \qquad 
\kappa (dx)= N(x, \{\partial\}) \mu_H (dx).
$$
Here $\mu_H$ is the Revuz measure of the PCAF $H$ and $\partial$ is the cemetery point for $X$
added to $E$ as a one-point compactification.  

 Let $F=\cup^{N}_{j=1}K_j$ be the union of $N$ disjoint compact subsets $K_j$ of positive $\EE$-capacity. Set
 $D:=E\setminus F$. In this paper, we will construct a new Markov process $X^*$ from $X$  by darning 
(or shorting) each hole $K_j$ into a single point $a^*_j$. This new process has state space $E^*:=D\cup\{a^*_1,\cdots, a^*_N\}$ and is $m^*$-symmetric, where $m^*:=m$ on $D$ and $m^*(E^* \setminus D):=0$.
Moreover, the jumping measure $J^*$ and the killing measure $\kappa^*$ of $X^*$ on $E^*$ 
should have the properties
inherited from $J$ and $\kappa$ without incurring additional jumps or killings; that is,
\begin{equation}\label{e:1.1}
J^*=J  \hbox{ on } D\times D, \quad J^*(a^*_i, dy)=J(K_i, dy)  \hbox{ on }  D, 
\quad J^*(a^*_i, a^*_j)=J(K_i, K_j)  \hbox{ for } i\not=j, 
\end{equation}
\begin{equation}\label{e:1.2}
\kappa^*=\kappa \quad\hbox{on } D \quad \hbox{ and } \quad \kappa^*(a^*_j)=\kappa (K_j)
\hbox{ for } 1\leq j\leq N.
\end{equation}
We will show that such $X^*$ always exists and is unique in law. 
This process $X^*$ coincides with the Markov process with  darning introduced in 
\cite{CF} under the assumption that $X$ enters each $K_j$ in a continuous way, that is, 
$X_{\tau_D-} \in F$ on $\{\tau_D <\zeta\}$; see \cite[Theorem 7.7.3]{CF}. 
Here $\zeta$ is the lifetime of $X$ and $\tau_D:=\inf\{t\geq 0: X_t\notin D\}$ is the first
exit time from $D$ by the process $X$. 
Thus we will call $X^*$ the Markov process obtained from $X$ by darning (or shorting)
each $K_j$ into a singleton $a_j^*$, or simply, Markov process with darning.   
Note that as a consequence of the $m^*$-symmetry assumption,
$X^*$ spends zero Lebesgue amount of time on $E^*\setminus D=\{a^*_1, \cdots, a^*_N\}$.

The (new) Markov process with darning $X^*$ will be constructed from $X$ via Dirichlet form method. 
Since in applications,  $\EE (u, u)$ can be interpreted as energy of a potential $u\in \FF$, 
intuitively speaking,  restricting $\EE$ to those $u\in \FF$ that
are constant $\EE$-q.e. on each $K_j$ exactly represents shorting each $K_j$ 
into a single point $a^*_j$.
Our Theorem \ref{T:3.3} of this paper shows that, after a suitable identification,
 this approach indeed works in great generality, 
without any additional assumptions.  We will further show in Theorem \ref{T:3.4}  that it is unique in distribution. 
 When $X$ is a Brownian motion in $\R^n$ and $F=K$ is a compact set, the above Dirichlet form
method of constructing $X^*$ was carried out in  \cite{Chen, CFR} and we call $X^*$ Brownian motion
with darning (BMD). 
When $E$ is the exterior of the  unit disk in $\R^2$, $F=\partial E$
  and $X$ is the reflecting Brownian motion on $E$, BMD $X^*$ has the same law as
  the excursion reflected Brownian motion appeared in \cite{La} in connection with the study of $SLE$ in multiply  connected planar domains. Planar BMD enjoys conformal invariance property, see \cite{CF}.
In \cite{CFR, CF2,  CF3},  BMD has been used to study Chordal Komatu-Loewner equation and 
stochastic Komatu-Loewner equation in standard slit domains in upper half space.

	The second goal of this paper is to present approximation schemes for general symmetric 
  Markov processes with darning $X^*$ in the finite dimensional sense, which can also be used
	to simulate the darning processes. 
We note that the construction of $X^*$ either by Dirichlet form method as in this paper 
or by Poisson point process of excursions when the process $X$ enters the holes in
a continuous way as in \cite{CF, CFY, FT} does not provide a practical way to simulate $X^*$.
 Our approach of this paper is to introduce additional jumps among each $K_j$ with large intensity.
 Intuitively, when the jumping intensity for these additional jumps increases to infinity, the new process can no longer
  distinguish points among each $K_j$,
 which would result in shorting (or darning) each $K_j$ into a single point $a^*_j$.
 To be precise, for each $j$, let $\mu_j$ be a finite smooth whose quasi-support is $K_j$ and having bounded 1-potential $G_1 \mu_j$. For each $\lambda >0$, consider the following Dirichlet form $(\mathcal{E}^{(\lambda)},\mathcal{F})$ on $L^2(E; m)$:
\begin{equation}\label{e:1.3}
 \mathcal{E}^{(\lambda)}(u,v)=\mathcal{E}(u,v)+\lambda \sum_{j=1}^N\int_{K_j\times K_j}(u(x)-u(y))(v(x)-v(y))
 \mu_j(dx)\mu_j(dy)
 \end{equation}
 for $u,v\in \mathcal{F}$. It is easy to see that $(\mathcal{E}^{(\lambda)},\mathcal{F})$ is a regular Dirichlet form on $L^2(E; m)$ and thus by \cite{FOT}, there is a $m$-symmetric Hunt process $X^{(\lambda)}$ associated with it.
The process $X^{(\lambda)}$ is the superposition of $X$ with jumps among points within each $K_j$.
 The process $X^{(\lambda)}$ can also be obtained from $X$ by the following piecing together procedure.
 Let $X^0$ be the subprocess obtained from $X$
 through killing via measure $\lambda \sum_{j=1}^N \mu_j$. More precisely, let $A_s^j$ be positive continuous additive functional (PCAF in abbreviation) of $X$ with Revuz measure $\mu_j$. Then the law of $X^0$ is determined by the following: for every positive function $f$ on $E$,
$$ 
\E_x \left[ f(X_t^0) \right]=\E_x \left[ e^{-\lambda \sum_{j=1}^NA_t^j}f(X_t) \right] .
$$
 Denote by $\zeta^0$ the lifetime of $X^0$. For each starting point $x\in E$, $X^{(\lambda)}$ can be obtained
 from $X^0$ through the following redistribution and patching procedure.
 Run a copy of $X^0$ starting from $x$ and set $X_t^{(\lambda)}=X_t^0$ for $t\in [0,T_1)$, where $T_1=\zeta_1^0$ is the lifetime of $X^0$ starting from $x$. If $\zeta^0=\infty$ or $X^{(\lambda)}_{T_1-}=\partial$, then we define $X_t^{(\lambda)}=\partial$ for $t\geq T_1$. Otherwise, $X_{T_1-}^{(\lambda)}\in F$, say
 $X_{T_1-}^{(\lambda)}\in K_{j_1}$. Select $x_2\in K_{j_1}$ according to the probability distribution $\mu_{j_1}/\mu_{j_1}(K_j)$ and define $X_{T_1-}^{(\lambda)}=x_2$. Run an independent copy of $X^0$ starting from $x_2$, whose lifetime is denoted as $\zeta_2^0$. Define $X_{T_1+t}^{(\lambda)}=X_t^0$ for $t\in [0,\zeta_2^0)$
 and set $T_2=T_1+\zeta_2^0$. If $T_2=\infty $ or $X_{T_2-}^{(\lambda)}=\partial$, then we define $X_t^{(\lambda)}=\partial$ for $t\geq T_2$. Otherwise, $X_{T_2-}^{(\lambda)}\in F$, say $X_{T_2-}^{(\lambda)}
 \in K_{j_2}$. Let $x_3\in K_{j_2}$ according to the probability distribution $\mu_{j_2}/\mu_{j_2}(K_{j_2})$ and define $X_{T_2}^{(\lambda)}=x_3$, and so on. The above described patching together procedure is a particular
 case discussed in \cite{INW}. The resulting process is a Hunt process on $E$. It is easy to verify that the Dirichlet
 form associated with $X^{(\lambda)}$ is $(\mathcal{E}^{(\lambda)}, \mathcal{F})$.

When the intensity $\lambda \to \infty$, process $X^{(\lambda)}$ behaves like $X$ outside $F$ but
can not distinguish points in each $K_j$. In other words, in the limit, each $K_j$ is collapsed  into a single point $a^*_j$.
So if the limit exists, the limiting process should be Markov process with darning of $X$ but up to a time change.
This is because $m$ is a symmetrizing measure for each $X^{(\lambda)}$ so under stationarity,
each $X^{(\lambda)}$ spends time in $F$ at a rate proportional to $m(F)$.
Let $Y$ be the Hunt process on $E^*$ obtained from $X^*$ through a time change via Revuz measure
$\mu = m|_D + \sum_{j=1}^N m(K_j) \delta_{\{ a^*_j\}}$. 
That is, $Y$ is a sticky Markov process with darning, which spends $m(K_j)$-proportional
 Lebesgue amount of time at $a^*_j$ over time duration.
 We show in this paper that as $\lambda\rightarrow \infty$, $X^{(\lambda)}$ converges to $Y$ 
in the finite dimensional sense;
 see Theorem \ref{T:4.3} below for a precise statement.

When $X$ is a diffusion process on $E$ and each compact set $K_j$ is connected, it is possible
to get the diffusion with darning $X^*$ on $E^*$, up to a time change, by increasing the   conductance 
on each $K_j$ to infinity. This is illustrated by Theorem \ref{T:4.4}.

An effective way of establishing finite dimensional convergence for symmetric Markov processes is the
Mosco convergence of Dirichlet forms \cite{Mo}. However, the state space $E^*$ of $X^*$ is different from
that of $X^{(\lambda)}$ -- there is a sudden collapsing of the state space 
right at the limit $\lambda = \infty$. 
This is in stark contrast with cases considered in \cite{CKK, Ki, KS}, 
where the weak converges and Mosco convergenc are studied for processes
and for Dirichlet forms on different state spaces.
In these papers, the state spaces are changing in a continuous way 
as $n\to \infty$. 
From the Dirichlet form point of view, the domain $\FF^*$ of the Dirichlet form $(\EE^*, 
\FF^*)$ associated with our $X^*$, 
viewed as a subspace of $L^2(E; m)$, is exactly those in $\FF$ that are constant quasi-everywhere
on each $K_j$. So it may  not be dense in $L^2(E; m)$ in general while the domain 
of the Dirichlet form for  $X^{(\lambda)}$ is $\FF$ for every $\lambda >0$.
Thus the existing theory of Mosco convergence \cite{Mo, KS, Ki, CKK} can not be applied directly.
In Section \ref{S:2}  of this paper, we extend the characterization of $L^2$-convergence of semigroups
for Mosco convergence to closed symmetric  forms whose domains may not be dense in $L^2$-space;
see Theorem \ref{T:2.3}. This result may be of independent interest. 
 The approximation schemes mentioned above for Markov processes with darning are established by
applying Theorem \ref{T:2.3}.

The idea and approach of this paper, including the generalized Mosco convergence result 
(Theorem \ref{T:2.3}), can also be used
to study approximation for other darning related processes. In Section \ref{S:5},
we illustrate how to use the ideas of this paper to approximate Brownian motion in a plane with a very thin flag pole studied recently in \cite{CL}
by Brownian motion in the plane with a vertical cylinder whose horizontal motion on the cylinder 
is a circular Brownian motion  moving at  fast speed.

\section{Mosco convergence of general  closed symmetric  forms}\label{S:2}

 One way to establish the finite dimensional convergence is via Mosco convergence \cite{Mo}.  However, the characterization of 
 convergence of symmetric semirgoups in \cite{Mo} is formulated
 only for those closed symmetric    forms whose domains of definition are dense in the $L^2$-spaces. In this section, we
  study Mosco convergence of general closed symmetric  forms  whose domains of definition may not necessarily be dense in the 
  corresponding $L^2$-spaces.

 Let $E$ be a locally compact separable metric space and $m$ a Radon measure on $E$ with full support.
 Suppose $(\EE, \FF)$ is a closed symmetric  form on $L^2(E; m)$; that is, $\FF$ is a linear subspace of $L^2(E; m)$,
 $\EE$ is a non-negative definite symmetric  form
 defined on $\FF\times \FF$ such that $\FF$ is a Hilbert space with inner product $\EE_1$. Here for $\alpha >0$,
 $$
 \EE_\alpha (f, g):= \EE (f, g) + \alpha (f, g)_{L^2(E; m)}, \qquad f, g \in \FF.
 $$
Note that here we do not assume $\FF$ is dense in $L^2(E; m)$.
 Throughout this paper, we use the convention that we define $\EE (f, f)=\infty$ for $f\notin \FF$.
  Given a    closed symmetric form $(\EE, \FF)$ on $L^2(E; m)$, by Riesz representation theorem, for every $f\in L^2(E; m)$ and $\alpha >0$,
  there is a unique $G_\alpha f \in \FF$ such that
 \begin{equation}\label{e:2.1a}
   \EE_\alpha (G_\alpha f , g) = (f, g)_{L^2(E; m)}
    \quad \hbox{for every } g\in \FF.
 \end{equation}
 These linear operators $\{G_\alpha, \alpha >0\}$ on $L^2(E; m)$ is called the resolvent of $(\EE, \FF)$.
 It is known that the resolvent $\{G_\alpha, \alpha>0\}$ of $(\EE, \FF)$ is strongly continuous
(that is, $\lim_{\alpha \to \infty} \| \alpha G_\alpha f -f \|_{L^2(E; m)}=0$ for every $f\in L^2(E; m)$)
 if and only if $\FF$ is dense in $L^2(E; m)$.
 If $\FF$ is dense in $L^2(E; m)$, then there is a unique, strongly continuous semigroup
 $\{P_t, t\geq 0\}$ associated with the strongly continuous resolvent $\{G_\alpha, \alpha >0\}$,
 and hence with $(\EE, \FF)$.

If $\FF$ is not dense in $L^2(E; m)$,
 denote by $\overline \FF$ the closure of $\FF$ in $L^2(E; m)$.
Then $(\EE, \FF)$ is a closed symmetric  form on $\overline \FF$.
The following facts are known; see \cite[pp.2-4]{CF} or \cite{MR}.
There is a unique strongly continuous contraction symmetric resolvent $\{\wh G_\alpha; \alpha >0 \}$
on $\overline \FF$ associated with it:
\begin{equation}\label{e:2.2a}
   \EE_\alpha (\wh G_\alpha f , g) = (f, g)_{L^2(E; m)}
    \quad \hbox{for every } g\in \FF.
 \end{equation}
It in turn is associated with
a unique strongly continuous contraction symmetric semigroup $\{\wh P_t; t\geq 0\}$
on $\overline \FF$ via
$$
\wh G_\alpha f = \int_0^\infty e^{-\alpha t} \wh P_t f dt, \quad f\in \overline \FF.
$$
The correspondence between $(\EE, \FF)$, $\{\wh G_\alpha, \alpha >0\}$
and $\{\wh \P_t, t\geq 0\}$  on $\overline \FF$ are one-to-one. In particular,
\begin{eqnarray*}
 \FF &=& \left \{u\in \overline \FF :\lim_{t\rightarrow 0}\frac{1}{t}(u- \wh P_tu,u)_{L^2(E; m)}<\infty\right\}, \\
\EE (u,v)&=& \lim_{t\rightarrow 0}\frac{1}{t}(u- \wh P_tu,v)_{L^2(E; m)}
\quad \hbox{for } u, v \in \FF.
\end{eqnarray*}
Denote by $(\wh \LL, \DD (\wh \LL))$ the generator of $\{\wh P_t; t\geq 0\}$ in
the Hilbert space $\overline \FF$ (equipped with
the $L^2$-inner product from $L^2(E; m)$). Then $u\in \DD (\wh \LL)$ if and only $u\in \FF$ and there is $f\in \overline \FF$
so that
$$ \EE (u, v) = -(f, v)_{L^2(E; m)} \quad \hbox{for every } v\in \FF;
$$
in this case, $\wh \LL u =f$. We have $G_\alpha (\overline \FF)=\DD (\LL)$ and $\wh P_t (\overline \FF) \subset \DD (\LL)$.

  Let $\Pi$ be the orthogonal projection operator from $L^2(E; m)$ onto $\overline \FF$.
Then we have from \eqref{e:2.1a} and \eqref{e:2.2a} that
\begin{equation}\label{e:2.3a}
G_\alpha f = \wh G_\alpha (\Pi f)  \quad \hbox{for every } \alpha >0 
\hbox{ and } f\in L^2(E; m).
\end{equation}

 \begin{defn}\label{d:2.1} A sequence of closed symmetric  forms $\{(\mathcal{E}^n,\mathcal{F}^n)\}$ on $L^2(E; m)$
 is said to be convergent to a closed  symmetric form $(\mathcal{E},\mathcal{F})$ on $L^2(E; m)$ in the sense of
 Mosco  if

 (a) For every sequence $\{u_n,n\geq 1\}$ in $L^2(E; m)$ that converges weakly to $u$ in $L^2(E;m)$,

   $$\liminf_{n\rightarrow \infty}\mathcal{E}^n(u_n,u_n)\geq \mathcal{E}(u,u),$$

 (b) For every $u\in L^2(E; m)$, there is a sequence $\{u_n,n\geq 1\}$ in $L^2(E; m)$ converging strongly to $u$
    in $L^2(E; m)$ such that

    $$\limsup_{n\rightarrow \infty} \mathcal{E}^n(u_n,u_n)\leq \mathcal{E}(u,u).$$
\end{defn}

 Denote by  $\{G_{\alpha}, \alpha>0\}$, $\{G_{\alpha}^n, \alpha>0\}$ the corresponding resolvents
 of $(\EE, \FF)$ and $(\EE^n, \FF^n)$, respectively.
 When $\FF$ and $\FF^n$ are dense in $L^2(E; m)$, the associated semigroup will be denoted by
 $\{P_t, t \geq 0\}$ and $\{P_t^n,t\geq 0\}$, respectively.

 The following result is known (see Theorem 2.4.1 and Corollary 2.6.1 of \cite{Mo}).

\begin{prop}\label{P:2.2} Let $(\mathcal{E},\mathcal{F})$ and $\{(\mathcal{E}^n,\mathcal{F}^n),n\geq 1)\}$
 be a sequence of closed symmetric    forms on $L^2(E; m)$. The following are equivalent:

\begin{description}
 \item{\rm (i)}  $(\mathcal{E}^n,\mathcal{F}^n)$ converges to $(\mathcal{E},\mathcal{F})$ in the sense of Mosco; 

 \item{\rm (ii)} For every $\alpha>0$ and $f\in L^2(E; m)$, $G_{\alpha}^n f$ converges to $G_{\alpha}f$ in $L^2(E; m)$
 as $n\to \infty$;

  \item{\rm (iii)}
  When $\FF^n$ and $\FF$ are all dense in $L^2(E; m)$, then (i) is equivalent to the following: For every $t>0$ and $f\in L^2(E; m)$, $P_t^nf$ converges to $P_tf$ in $L^2(E; m)$ as $n\to \infty$.
\end{description}
\end{prop}

The next result addresses the case when $\FF^n$ and  $\FF$ may not be dense in $L^2(E; m)$.

\begin{thm}\label{T:2.3}
 Let $(\mathcal{E},\mathcal{F})$ and $\{(\mathcal{E}^n,\mathcal{F}^n),n\geq 1)\}$
 be  closed symmetric  forms on $L^2(E; m)$. Let $\overline \FF^n$ and $\overline \FF$ be the closure of
 $\FF^n$ and $\FF$ in $L^2(E; m)$, respectively. Suppose that $\overline \FF^n \supset \overline \FF$ for every $n\geq 1$.
 Let $(\wh P^n_t; t\geq 0\}$ and $(\wh P_t; t\geq 0\}$ be the semigroups on $\overline \FF^n$ and $\overline \FF$
 associated with $(\EE^n, \FF^n)$ and $(\EE, \FF)$, respectively.
 Then the following hold.

\begin{description}
 \item{\rm (i)}  If $(\mathcal{E}^n,\mathcal{F}^n)$ converges to $(\mathcal{E},\mathcal{F})$ in the sense of Mosco,
 then for every $t>0$ and $f\in \overline \FF$, $\wh P_t^nf$ converges to $\wh P_tf$ in $L^2(E; m)$ as $n\to \infty$.

  \item{\rm (ii)} Suppose that the closed subspace $\overline \FF^n$ converges to $\overline \FF$
	in $L^2(E; m)$ in the sense that
$$
\lim_{n\to \infty} \| \Pi^n f -\Pi f\|_{L^2(E; m)}=0
\quad \hbox{for every } f\in L^2(E; m),
$$
where $\Pi^n$ and $\Pi$ denote the orthogonal projection operators of $L^2(E; m)$ onto $\overline \FF^n$ and $\overline \FF$,
respectively.
If  $\wh P_t^nf$ converges to $\wh P_tf$ in $L^2(E; m)$ for every $t>0$ and $f\in \overline \FF$,
  then $(\mathcal{E}^n,\mathcal{F}^n)$ converges to $(\mathcal{E},\mathcal{F})$ in the sense of Mosco.
\end{description}
\end{thm}

\pf Let $\{G^n_\alpha; \alpha>0\}$ and $(\wh G^n_\alpha; \alpha >0\}$ be the resolvents
 on $L^2(E; m)$ and on $\overline \FF^n$, respectively,  associated with the closed symmetric form 
 $(\EE^n, \FF^n)$ via \eqref{e:2.1a} and \eqref{e:2.2a}.
Similar notations   $\{G_\alpha; \alpha>0\}$ and $(\wh G_\alpha; \alpha >0\}$ will be used for $(\EE, \FF)$.
We know from Proposition \ref{P:2.2} that  $(\mathcal{E}^n,\mathcal{F}^n)$ converges to $(\mathcal{E},\mathcal{F})$ in the sense of Mosco if and only if $G^n_\alpha f$ converges to $G_\alpha f$
in $L^2(E; m)$ for every $f\in L^2(E; m)$.

(i) Suppose that $(\mathcal{E}^n,\mathcal{F}^n)$ converges to $(\mathcal{E},\mathcal{F})$ in the sense of Mosco,.
Then in view of \eqref{e:2.3a} and the assumption that $\overline \FF^n \supset \overline \FF$,
we have for every $\alpha >0$ and $f\in \overline \FF$, $\wh G^n_\alpha f$ converges to $\wh G_\alpha f$
in $L^2(E; m)$. We claim this implies that $\wh P_t^nf$ converges to $\wh P_tf$ in $L^2(E; m)$ as $n\to \infty$
for every $t>0$ and $f\in L^2(E; m)$. The proof is  similar to that for \cite[Theorem IX.2.16]{Ka}.
For reader's convenience, we spell out the details here.

Denote by $\wh \LL^n$ and $\wh \LL$ the generators of the strongly continuous semigroups $\{\wh P^n_t; t\geq 0\}$ and
$\{\wh P_t; t\geq 0\}$,  respectively.  Note that
\begin{eqnarray*}
 \frac{d}{dt} \wh P_t \wh G_\alpha
 &=& \wh P_t \wh \LL \wh G_\alpha = \wh P_t (\alpha \wh G_\alpha -I), \\
 \frac{d}{dt} \wh P^n_t \wh G^n_\alpha
 &=& \wh P^n_t \wh \LL^n \wh G^n_\alpha = \wh P^n_t (\alpha \wh G^n_\alpha -I).
\end{eqnarray*}
Thus in view of $\overline \FF \subset \overline \FF^n$, we have
\begin{eqnarray*}
\frac{d}{ds} \wh P^n_{t-s} \wh G^n_\alpha \wh P_s \wh G_\alpha
&=& - \wh P^n_{t-s} (\alpha \wh G^n_\alpha -I) \wh P_s \wh G_\alpha
   +  \wh P^n_{t-s} \wh G^n_\alpha \wh P_s (\alpha \wh G_\alpha -I) \\
&=&  \wh P^n_{t-s} ( \wh P_s \wh G_\alpha -  \wh G^n_\alpha \wh P_s )
=  \wh P^n_{t-s} ( \wh G_\alpha -  \wh G^n_\alpha ) \wh P_s .
\end{eqnarray*}
Integrating in $s$ over $[0, t]$ yields
$$
 \wh G^n_\alpha \wh P_t \wh G_\alpha
-\wh P^n_{t} \wh G^n_\alpha  \wh G_\alpha
= \int_0^t \wh P^n_{t-s} ( \wh G_\alpha -  \wh G^n_\alpha ) \wh P_s  ds . 
$$
Hence for every $f\in \overline \FF$,
$$
\lim_{n\to \infty}
\| \wh G^n_\alpha (\wh P_t -\wh P^n_t) \wh G_\alpha f \|_{L^2(E;m)}
\leq \lim_{n\to \infty} \int_0^t \| \wh P^n_{t-s} ( \wh G_\alpha -  \wh G^n_\alpha )
 \wh P_s f\|_{L^2(E; m)} ds =0.
$$
Since $\wh G_\alpha (\overline \FF)$ is $L^2$-dense in $\overline \FF$, we have for every $u\in \overline \FF$,
$$
\lim_{n\to \infty}
\| \wh G^n_\alpha (\wh P_t -\wh P^n_t) u \|_{L^2(E;m)} =0.
$$
On the other hand,
by the $L^2$-contraction property of $\wh P^n_t$ and $\wh P_t$, we have
 for $u\in \overline \FF \subset \overline \FF^n$,
$ \wh G^n_\alpha \wh P_t u -\wh P_t \wh G_\alpha u
=(\wh G^n_\alpha -\wh G_\alpha ) \wh P_t u \to 0$ in $L^2(E; m)$ as $n\to \infty$,
and $ \wh G^n_\alpha \wh P^n_t u -\wh P^n_t \wh G_\alpha u
=\wh P^n_t (\wh G^n_\alpha -\wh G_\alpha ) u\to 0$ in $L^2(E; m)$ as $n\to \infty$.
It follows then
$$ \lim_{n\to \infty} \| (\wh P^n_t -\wh P_t) \wh G_\alpha u\|_{L^2(E; m)} =0
\quad \hbox{for every } u\in \overline \FF.
$$
Since $\wh G_\alpha (\overline \FF)$ is $L^2$-dense in $\overline \FF$, we have
$ \lim_{n\to \infty}  \| (\wh P^n_t -\wh P_t) u \|_{L^2(E; m)} =0$ for every
$u\in \overline \FF$.

(ii) Conversely, assume $\overline \FF^n$  converges to  $ \overline \FF$
and $ \lim_{n\to \infty}  \| (\wh P^n_t -\wh P_t) u \|_{L^2(E; m)} =0$ for every
$u\in \overline \FF$. Denote by $\Pi^n$ and $\Pi$ the orthogonal projection operator of $L^2(E; m)$ onto
$\overline \FF^n$ and $\overline \FF$, respectively.
We have by \eqref{e:2.3a} and the $L^2$-contraction property of $\wh \P^n_t$ and $\wh P_t$ that for
every $\alpha>0$ and $f\in L^2(E; m)$,
\begin{eqnarray*}
\lim_{n\to \infty} \|G^n_\alpha f - G_\alpha f\|_{L^2(E; m)}
&=& \lim_{n\to \infty} \|\wh G^n_\alpha (\Pi^n f) - \wh G_\alpha (\Pi f)\|_{L^2(E; m)} \\
&\leq&  \lim_{n\to \infty} \left( \| ( \wh G^n_\alpha - \wh G_\alpha) (\Pi f ) \|_{L^2(E; m)} +
 \| \wh G^n_\alpha (\Pi^n f -\Pi f)\|_{L^2(E; m)} \right) \\
&\leq & \lim_{n\to \infty} \alpha^{-1} \| \Pi^n f -\Pi f\|_{L^2(E; m)} =0 .
\end{eqnarray*}
It follows from Proposition \ref{P:2.2}
that $(\mathcal{E}^n,\mathcal{F}^n)$ converges to $(\mathcal{E},\mathcal{F})$ in the sense of Mosco.
\qed

\section{Markov processes with darning }\label{S:3}

 Suppose $(\EE, \FF)$ is a regular symmetric Dirichlet form on $L^2(E; m)$.
 In particular, $\FF$ is a dense linear subspace of $L^2(E; m)$.
Let $X$ be the Hunt process on $E$ associated with $(\EE, \FF)$. 
In the remainder of this paper, we use the convention that every $f\in \FF$ is represented by its quasi-continuous version,
 which is unique up to an $\EE$-polar set.
Suppose that  $K_1, \dots, K_N$ are separated, non-$\EE$-polar compact subsets of $E$.
Let $F=\cup_{j=1}^N K_j$ and $D=E\setminus F$. We short (or collapse) each $K_j$ 
into a single point $a^*_j$.
Formally, by identifying each $K_j$ with a single point $a^*_j$, we
can get an induced topological space $E^*:=D\cup \{a_1^*, \dots,
a^*_N\}$ from $E$, with a neighborhood of each $a^*_j$ defined
as $(U\cap D)\cup \{a^*_j\}$ for some neighborhood $U$ of
$K_j$ in $E$. Let $m^*=m$ and $D$ and set $m^* (F^*)=0$, where $F^*:=\{a_1^*, \dots, a_N^*\}$.

\begin{defn}\label{D:3.1} \rm  A strong Markov process  on $E^*$
is said to be Markov process with darning obtained from $X$ by shorting each $K_j$ into a single point $a^*_j$, or simply a Markov process with darning,  is an $m^*$-symmetric Markov process $X^*$ on $E^*$ such that
\begin{description}
\item{(i)} the part process  of $X^*$ in $D$ has the same law as the part process $X$ in $D$ for 
$\EE$-q.e. starting point in $D$; 
\item{(ii)} The jumping measure $J^*(dx, dy)$ and killing measure $\kappa^*$ of $X^*$ on $E^*$ 
have the properties inherited from $X$ without incurring additional jumps or killings, that is, they have the properties  
\eqref{e:1.1} and \eqref{e:1.2}. 
\end{description}
\end{defn}

\begin{remark}\label{R:3.2}  \rm
Note that if $X^*$ is a Markov process with darning of $X$,
it follows from Definition \ref{D:3.1} that 
\begin{equation}\label{e:3.1} 
\P_x (X^*_{\tau_D}= a^*_j)= \P_x (X_{\tau_D}\in K_j)
\quad \hbox{for q.e. } x\in D,
\end{equation} 
where $\tau_D:=\inf\{t>0: X^*_t\notin D\}=\inf\{t>0: X_t\notin D\}$. 
Hence each $a^*_j$ is of positive capacity with respect to the process $X^*$ because 
$K_j$ is of positive $\EE$-capacity. 
In particular,  each $a^*_j$  is regular for itself;  that is, $\P_{a^*_j}(\sigma_{a^*_j}=0 )=1$,
where $\sigma_{a^*_j}:=\inf\{t>0: X^*_t=a^*_j\}$.  
This is due to the general fact that for any nearly Borel measurable set  $A\subset E^*$, 
$A \setminus  A^r$ is semipolar for process $X^*$ and hence $\EE^*$-polar. 
Here $A^r$ denotes all the regular points for $A$ with respect to the strong Markov process $X^*$. 
\end{remark}

\medskip

We will show in this section that Markov process with darning $X^*$ from $X$ always exists and is unique in distribution.

\medskip

   For $1\leq j\leq N$ and $\alpha>0$, define
 $$
  \varphi^{(j)}(x):=\P_x(X_{\sigma_F}\in K_j) \quad \hbox{ and } \quad
  u_{\alpha}^{(j)}(x):=\E_x[e^{-\alpha \sigma_F};X_{\sigma_F}\in K_j].
 $$
 Here $\sigma_F:=\inf\{t\geq 0:X_t\in F\}$. Let $\mathcal{H}_{\alpha}$ be the linear span of $\{u_{\alpha}^{(j)}, j=1,\cdots,N\}$, and  $(\mathcal{E},\mathcal{F}^D)$  the Dirichlet form for the  part process $X^D$
of $X$ killed upon exiting $D$.
Since each $K_j$ is compact and $(\mathcal{E},\mathcal{F})$ is a regular Dirichlet form, there is a function $f_j\in C_c(E)\cap \mathcal{F}$ such that $f_j=1$ on $K_j$ and
 $f_j=0$ on $\cup_{l:l\neq j}K_l$.
Consequently, $u_{\alpha}^{(j)}(x)=\E_x[f_j(X_{\sigma_K})]$ is the $\mathcal{E}_{\alpha}$-orthogonal projection of $f_j$  to the complement of $\mathcal{F}^D$, where
$$
\EE_\alpha (u, v):= \EE (u, v) + \alpha (u, v)_{L^2(E; m)}
\quad \hbox{for } u, v \in \FF.
$$
 So in particular, $u_{\alpha}^{(j)} \in \FF$ for every $\alpha >0$ and $1\leq j \leq N$.   Define
\begin{equation}\label{e:3.2}
 \mathcal{F}^*=\mathcal{F}^D\oplus \mathcal{H}_{\alpha} .
 \end{equation}
It is easy to see that the above definition of $\FF^*$ is independent of $\alpha>0$.
The space
 $\mathcal{F}^D$ is exactly the collection of functions in $\mathcal{F}$ that vanish $\EE$-quasi-everywhere 
($\EE$-q.e. in abbreviation) on $D^c$,
while $u^{(j)}_\alpha =1$ on $\EE$-q.e. on $K_j$ and vanishes $\EE$-q.e. on $K_l$ for $l\not= j$. 
Hence by regarding each $u^{(j)}_\alpha$ as a function defined on $E^*$, 
$\FF^*$ can be viewed as a dense linear subspace of $L^2(E^*; m^*)$.
Define 
\begin{equation}\label{e:3.3}
\EE^*(u, v)= \EE (u, v) \quad \hbox{for } u, v \in \FF^*.
\end{equation}
We will show in Theorem \ref{T:3.3} below that
$(\EE^*, \FF^*)$ is a regular Dirichlet form on $L^2(E^* ;  m^*)$.
Consequently, it uniquely determines a Hunt process $X^*$ on $E^*$.

As we saw from the above, 
 $\mathcal{F}^*=\mathcal{F}^D\oplus \mathcal{H}_{\alpha}$ can be identified with 
functions in $\mathcal{F}$ that are constant $\EE$-q.e. on each $K_j$.
For $f\in \FF$, define $\H_F^1 f(x)=\E_x \left[ e^{-\sigma_F} f(X_{\sigma_F})\right]$. 
Note that $f-\H^1_F f \in \FF^D$ and
$\H^1_F f$ is $\EE_1$-orthogonal to $\FF^D$.  

\begin{thm}\label{T:3.3}
 $(\EE^*, \FF^*)$ is a regular symmetric Dirichlet form on $L^2(E^*; m^*)$
and its associated Hunt process $X^*$ on $E^*$ 
is a Markov process of darning obtained from $X$ by shorting
each $K_j$ into a single point $a^*_j$.
\end{thm}

\pf  Let $\cC=\{u\in \FF \cap C_c(E): u \hbox{ is constant
on each } K_j\}$. By defining $u(a^*_j)$ to be the value
of $u$ on $K_j$, we can view $\cC$ as a subspace of
$\FF^* \cap C_c(E^*)$. Since $\cC$ is an algebra
that separates points in $E^*$, 
$\cC$ is uniformly dense in $C_\infty (E^*)$ by Stone-Weierstrass theorem.
Next we show that $\cC$ is $\EE_1^*$-dense in $\FF^*$.
For this, it suffices to establish that each $u^{(1)}_j(x):=\E_x \left[
e^{-\sigma_F}; X_{\sigma_F} \in K_j \right]$ can
be $\EE_1$-approximated by elements in $\cC$.
Let $f_j\in \FF\cap C_c(E)$ so that $f_j=1$ on $K_j$
and $f_j=0$ on $K_i$ for $i\not=j$. Note that
$u^{(1)}_j =\H_F^1 f_j = f_j-(f_j-\H^1_F f_j)$
and $f_j - \H^1_F f_j \in \FF^D$.
Since $(\EE, \FF^D)$ is a regular Dirichlet form
on $L^2(D; m)$, there is a sequence
$\{g_k, k\geq 1\}\subset \FF^D \cap C_c(D)$ that is $\EE_1$-convergent
to $f_j-\H_F^1 f_j$. Let $v_k:=f_j-g_k$, which is
in $\cC$ and $\EE^*_1$-convergent to $u^{(1)}_j$.
Thus we have established that $(\EE^*, \FF^*)$ is
a regular Dirichlet form on $L^2(E^*; m^*)$.

Let $X^*$ be the symmetric Hunt process on $E^*$
associated with the regular Dirichlet form $(\EE^*, \FF^*)$ on $L^2(E^*; m^*)$.
Clearly the part process of $X^{*, D}$ of $X^*$ in $D$ has the same distribution
as the part process of $X$ in $D$ because the part Dirichlet forms of $(\EE^*, \FF^*)$ and
$(\EE, \FF)$ on $D$ are the same. Denote by $J^*(dx, dy)$ and $\kappa^*$ the jumping measure
and killing measure of $(\EE^*, \FF^*)$. 
For every $f\in \FF^*$, by the Beurling-Deny decomposition of $(\EE^*, \FF^*)$, 
$$
\EE^*(f, f)= \EE^{* c} (f, f)+\frac12 \int_{E^*\times E^*} (f(x)-f(y))^2 J^*(dx, dy) 
+ \int_{E^*} f(x)^2 \kappa^* (dx),
$$
where $\EE^{* c}$ is a  non-negative definite symmetric bilinear form on $\FF^*$
that has strong local property.
On the other hand, by \eqref{e:3.2}, each $f\in \FF^*$ can be regarded  
a function in $\FF$ that is constant on each $K_j$ and 
\begin{eqnarray*}
\EE^*(f, f)=\EE (f, f)=\EE^c(f, f)+\frac12 \int_{E\times E} (f(x)-f(y))^2 J(dx, dy) 
 + \int_E f(x)^2 \kappa (dx).
\end{eqnarray*}
Comparing the above two displays yields $ \EE^{* c} (f, f) = \EE^c(f, f)$
and $J^*$ and $\kappa^*$ satisfy \eqref{e:1.1}-\eqref{e:1.2}. 
This proves that $X^*$ is a Markov process with darning for $X$.  \qed

The next result gives the uniqueness of the Markov process with darning for $X$.

\begin{thm}\label{T:3.4}
 Suppose $X^*$ is a Markov process with darning for $X$ in the sense of 
Definition \ref{D:3.1}. Then the Dirichlet form for $X^*$ on $L^2(E^*; m^*)$ 
is the one $(\EE^*, \FF^*)$ given by \eqref{e:3.2}-\eqref{e:3.3}. 
Consequently, Markov process with darning for $X$ is unique in distribution. 
\end{thm}

\pf Let $(\wt \EE, \wt \FF)$ be the quasi-regular Dirichlet form of $X^*$ on $L^2(E^*; m^*)$
(cf. \cite{CF, Fi}).  It suffices to show $(\wt \EE, \wt \FF) = (\EE^*, \FF^*)$.
 By Definition \ref{D:3.1}(i), $(\wt \EE, \wt \FF_D)=(\EE, \FF_D)$, 
where $\wt \FF_D$ and $\FF_D$ denote the part Dirichlet form of $(\wt \EE, \wt \FF)$
and $(\EE, \FF)$ on $D$, respectively.  
Let $F^*:=\{a^*_1,  \dots, a^*_N \}$ and
$\sigma^*:=\inf\{t>0: X^*_t\in F^*\}$.
By the $\EE^*_1$-orthogonal projection (see \cite[Theorem 3.2.2]{CF}), for every $u\in \FF^*$,
$\H^1_{F^*} u(x):=\E_x \left[ e^{-\sigma^*}
u(X^*_{\sigma^*}) \right] \in \FF$
and $u-\H^1_{F^*}u\in \wt \FF_D=\FF_D$.
It follows from Definition \ref{D:3.1} (cf. \eqref{e:3.1}) that,  for $x\in D$, 
$$
\H^1_{F^*}u(x)=
 \sum_{j=1}^N u(a_j^*)\, \E_x \left[ e^{-\sigma^*}; X^*_{\sigma^*}=a^*_j \right]
 =  \sum_{j=1}^N  u(a_j^*)\,  \E_x \left[ e^{-\sigma_F}; X_{\sigma_F}\in K_j \right]
 =  \sum_{j=1}^N u(a_j^*)\,  u_1^{(j)}(x) . 
$$
As  by  Remark \ref{R:3.2}, each $a^*_j$ is of positive $\EE^*$-capacity,  we have 
 $$
\{(u(a^*_1),  \dots, u (a^*_N)); u\in \FF\}=\R^N
$$
 and so  $\wt \FF=\FF^*$ by \eqref{e:3.1}.
 For $u\in \wt \FF =\FF^*$, 
 let $\wt \mu^{c}_{\<u\>}$ and $ \mu^{c}_{\<u\>}$ be the energy measure of $u$
 corresponding to the strongly local part $\wt \EE^c$ and $\EE^{* c}$ of the corresponding
 Dirichlet forms  $(\wt \EE, \wt \FF)$ and $(\EE^*, \FF^*)$, respectively.  
 Since $(\wt \EE, \wt \FF_D)=(\EE, \FF_D)=(\EE^*, \FF^*_D)$, we have
 $$
 \wt \mu^{c}_{\<u\>}(D)= \mu^{*c}_{\<u\>} (D).
 $$
 On the other hand,  for every bounded $u\in \wt \FF$,
 since the energy measures $\wt \mu^c_{\<u\>}$ and $\mu^{*c}_{\<u\>}$ of $u$
 do  not charge on level sets of $u$ (cf. \cite[Theorem 4.3.8]{CF}),  
 $$
 \wt \mu^{c}_{\<u\>} (a^*_j) = 0
 = \mu^{*c}_{\<u\>} (a^*_j) \quad \hbox{ for every } 1\leq j\leq N.
 $$
  Consequently, 
 $$ 
 \wt \mu^{c}_{\<u\>}(E^*\setminus D)=\sum_{j=1}^N \wt \mu^{c}_{\<u\>} (a^*_j)
 =0 = \sum_{j=1}^N   \mu^{*c}_{\<u\>} (a^*_j) = \mu^{*c}_{\<u\>}(E^*\setminus D) . 
 $$
 We conclude from the above two displays that 
 $$ 
 \wt \EE^c (u, u)=\frac12 \wt \mu^{c}_{\<u\>}(E^*) =  \frac12   \mu^{*c}_{\<u\>}(E^*) =\EE^{*c}(u, u) 
 $$
 for every bounded $u\in \wt \FF$ and hence for every $u\in \wt \FF$. 
 By the Beurling-Deny decomposition of $(\wt \EE, \wt \FF)$ and that $X^*$ is a Markov process 
 with darning for $X$,  we have for every $u\in \wt \FF=\FF^*$, 
 \begin{eqnarray*}
 \wt \EE(u, u)= \wt \EE^c (u, u) +  \frac12 \int_{E^*\times E^*} (f(x)-f(y))^2 J^*(dx, dy) 
+ \int_{E^*} f(x)^2 \kappa^* (dx) = \EE^* (u, u). 
 \end{eqnarray*} 
  This proves that  $(\wt \EE,  \wt \FF)=(\EE^*, \FF^*)$.
\qed 

\section{Approximation of Markov processes with darning}  \label{S:4}

We continue to work under the setting of Section \ref{S:3}. 
 Let  $X^*$ be the Markov process with darning obtained
from $X$ by shorting (or darning) each $K_j$ into a single point $a^*_j$. 
In this section, we study its approximations, whose scheme can be used to simulate $X^*$. 
For this, we first need to introduce sticky Markov process
with darning obtained from $X^*$ by a time change to possibly prolong the time spent on each $a^*_j$.

Define
$\mu= m^* + \sum_{j=1}^N m(K_j) \delta_{a^*_j}$, where $\delta_{a^*}$ is the Dirac measure concentrated at the point $a^*_j$. The smooth measure $\mu$ determines a positive continuous
additive functional $A^\mu$ of $X^*$. In fact, 
$$
A^\mu_t =t + \sum_{j=1}^N m(K_j) L^{a^*_j}_t,
$$ 
where $L^{a^*_j}$ is the local time of $X^*$ at $a^*_j$ having Revuz measure $\delta_{a^*_j}$.  
Let $\tau_t:=\inf\{s>0: A^\mu_s >t\}$ and $Y_t=X^*_{\tau_t}$.
Then the time-changed process $Y$ is $\mu$-symmetric and has Dirichlet form $(\EE^*, \FF^*)$ on $L^2(E^*; \mu)$;
see \cite{CF, FOT}.
The process $Y$ is a sticky Markov process with darning, as it may spend positive amount of Lebesgue time
at each $a^*_j$.

Conversely, starting with a sticky Markov process with darning $Y$ on $E^*$ associated with the regular Dirichlet form $(\EE^*, \FF^*)$ on $L^2(E^*; \mu)$, one can recover in distribution the Markov 
process with darning $X^*$ on $E^*$ through a time change as follows.
Let $\wt A_t=\int_0^t 1_D (Y_s)ds$, which is a positive continuous additive functional of $Y$
having Revuz measure $m^*$. Define its inverse $\wt \tau_t=\inf\{s>0: \wt A_s >t\}$.
Then $\wt X_t:=Y_{\wt \tau_t}$ is an $m^*$-symmetric strong Markov process on $E^*$ whose 
associated Dirichlet form is $(\EE^*, \FF^*)$ on $L^2(E; m^*)$ (cf. \cite{CF, FOT}).
In other words, $\wt X$ has the same distribution as $X^*$.

Let
\begin{equation}\label{e:4.1}
\wt \FF= \{ f\in \FF: f =\hbox{constant $\EE$-q.e. on } K_j \hbox{ for each } 1\leq j\leq N\} .
\end{equation}
Note that $(\EE, \wt \FF)$ is a closed symmetric  Markovian bilinear form on $L^2(E; m)$ but $\wt \FF$ is not
dense in $L^2(E; m)$ in general since each $K_j$ has positive $\EE$-capacity.
To emphasize its dependence on the domain of definition, we write $(\wt \EE, \wt \FF)$ for $(\EE, \wt \FF)$.
Denote by $\Pi$ the orthogonal projection of $L^2(E; m)$ onto the closure $\overline {\wt \FF}$ of $\wt \FF$ in $L^2(E; m)$.
Let $\{\wt G_\alpha, \alpha>0\}$ be the resolvent associated with
$(\wt \EE, \wt \FF)$ on $L^2(E; m)$, and $\{\wh P_t; t\geq 0\}$ and $(\wh G_\alpha; \alpha >0\}$
the semigroup and resolvent of the   closed symmetric  form $(\wt \EE, \wt \FF)$ on $\overline {\wt \FF}$, respectively.
We know from \eqref{e:2.3a} that $\wt G_\alpha f = \wh G_\alpha (\Pi f)$.
We now identify $\wh P_t$ and $\wh G_\alpha$, as well as $\Pi$.

The following map $T$ establishes a one-to-one and onto correspondence between the closed symmetric form
$(\wt \EE, \wt \FF)$ on $\overline {\wt \FF}\subset L^2(E; m)$
and the regular Dirichlet form $(\EE^*, \FF^*)$ on $L^2(E^*; \mu)$: for every $f\in \wt \FF$,
\begin{equation}
Tf=f \ \hbox{ on } D \quad \hbox{ and } \quad Tf(a^*_j)=f(K_j) \ \hbox{ for } 1\leq j\leq N.
\end{equation}
(For $g\in \FF^*$, $T^{-1} g(x)=g(x)$ for $x\in D$ and $T^{-1} g(x)= g(a^*_j)$  for $x\in K_j$.)
The map $T$ has the property that for every $f\in \wt \FF$,
\begin{equation}\label{e:4.3}
 \wt \EE (f, f)= \EE^*(Tf, Tf) \quad \hbox{and} \quad \| f\|_{L^2(E; m)}=\| Tf \|_{L^2(E^*; \mu)}.
\end{equation}
In other words, $T$ is an isometry between $(\wt \EE, \wt \FF)$ on $\overline {\wt \FF}\subset L^2(E; m)$
 and $(\EE^*, \FF^*)$ on $L^2(E^*; \mu)$
both in $\EE$ and in the $L^2$ sense.
Denote by $\{ G^*_\alpha; \alpha >0\}$ and $\{P^*_t; t\geq 0\}$ the resolvent and semigroup
associated with the regular Dirichlet form
$( \EE^*, \FF^*)$ on $L^2(E^*; \mu)$, and $\{ \wh G_\alpha; \alpha >0\}$ and $\{\wh P_t; t\geq 0\}$ the resolvent and semigroup on $\overline {\wt \FF}$
associated with the Dirichlet form
$(\wt \EE, \wt \FF)$ on $L^2(E; \mu)$.

For every $f\in L^2(E; m)$, we can define a function $f^*$ on $E^*$ by setting 
$f^*=f$ on $D$ and
\begin{equation}\label{e:4.4}
f^* (a^*_j) := \begin{cases}
\int_{K_j} f(y) m(dy)/m(K_j) & \hbox{when } m(K_j)>0, \cr
0  &\hbox{when } m(K_j)=0.
\end{cases}
\end{equation}
For an $\EE^*$-quasi-continuous function $g$  on $E^*$, we define
$$
T^{-1} g (x):= \begin{cases}
g(x) &\hbox{for } x\in D_0, \cr
g(a^*_j) & \hbox{for } x\in K_j .
\end{cases}
$$
Clearly, 
$$
T\circ T^{-1} g =g \qquad \hbox{and} \qquad 
\| T^{-1} g\|_{L^2(E; m)} =  \| g \|_{L^2(E^*; \mu)} . 
$$
Since $T$ extends to be an isometry between $\overline{\wt \FF}\subset L^2(E; m)$ and $L^2(E^*; \mu)$, 
the above defined $T^{-1}$ extends to be an isometry between $L^2(E^*; \mu)$ and $\overline{\wt \FF}\subset L^2(E; m)$.
We conclude that
$$
\overline {\wt \FF}=T^{-1}(L^2(E^*; \mu)) =\{f\in L^2(E; m): f \hbox{ is constant $m$-a.e. on each } K_i\}.
$$

\begin{thm}\label{T:4.1} \begin{description}
\item{\rm (i)}
For $f\in L^2(E; m)$,    $\Pi f = T^{-1} f^*$ $m$-a.e.

\item{\rm (ii)}  $\wt G_\alpha f =T^{-1} G^*_\alpha f^*$ for $f\in L^2(E; m)$.
\end{description}

\item{\rm (iii)} For $f\in \overline {\wt \FF}$,  $t>0$ and $\alpha>0$,
\begin{equation}\label{e:4.5}
\wh P_t f=  T^{-1} P^*_t (f^*) \quad \hbox{and} \quad  \wh G_\alpha f =
 T^{-1} G^*_\alpha (f^*).
\end{equation}
\end{thm}

\pf (i) Let ${\mathcal C}$ be the set defined in the proof of Theorem \ref{T:3.3}, which has shown to be
$\EE^*_1$-dense in $\FF^*$ in $L^2(E^*; \mu)$. So  in particular ${\mathcal C}$ is $L^2$-dense in $L^2(E^*; \mu)$.
Consequently, $T^{-1}{\mathcal C}$ is $L^2(E; m)$-dense in $\overline {\wt \FF}$.
 On the other hand, it is clear that for $f\in L^2(E; m)$,
$$
(f, T^{-1}g)_{L^2(E; m)} = (T^{-1} f^*, T^{-1}g)_{L^2(E; m)}
\quad \hbox{for every } g\in {\mathcal C}.
$$
Thus we have $\Pi f =T^{-1} f^*$ $m$-a.e.

(ii) Let $f\in L^2(E; m)$.
For every $\alpha >0$ and $g\in \wt \FF$, it follows from \eqref{e:4.3} that
\begin{eqnarray*}
\wt \EE_\alpha ( \wt G_\alpha f, g)
&=& \int_E f(x) g(x) m(dx) = \int_{E^*} f^*(x) T g(x) \mu (dx) \\
&=& \EE^* ( G^*_\alpha f^*, Tg)+\alpha \int_{E^*} G^* f^* (x) Tg (x) \mu(dx) \\
&=& \wt  \EE_\alpha (T^{-1} G^*_\alpha f^*, g).
\end{eqnarray*}
We thus conclude that $\wt G_\alpha f =T^{-1} G^*_\alpha f^*$.

(iii) This follows immediately from (i), (ii) and \eqref{e:2.3a} that for $f\in \overline {\wt \FF}$,
$$
\wh G_\alpha f =\wt G_\alpha f = T^{-1} G^*_\alpha f^*.
$$
It is clear that $T_t f:= T^{-1}P^*_t f^*$ defines a symmetric strongly continuous contraction semigroup
on $\overline {\wt \FF}\subset L^2(E; m)$, as
 $\{P^*_t; t\geq 0\}$ is a strongly continuous contraction semigroup on $L^2(E^*; m^*)$.
 Moreover, for every $\alpha>0$, $\int_0^\infty e^{-\alpha t} T_t f dt = T^{-1} G^*_\alpha f^*=\wh G_\alpha f$.
 Thus $T_t =\wh P_t$.
\qed

We now study an  approximation scheme of Markov processes with darning by introducing additional jumps over each $K_j$ with large intensity.
For each $j$, let $\mu_j$ be a finite smooth measure whose quasi-support is $K_j$ and having bounded 1-potential $G_1 \mu_j$, 
which always exists. For $\lambda >0$, consider the symmetric  regular Dirichlet form $(\EE^{(\lambda)}, \FF)$
  on $L^2(E; m)$ defined by \eqref{e:1.3}.
Observe that by \cite{SV} for every $u\in \FF$,
\begin{eqnarray*}
 \sum_{j=1}^N \int_{K_j\times K_j} (u(x)-u(y))^2 \mu_j (dx) \mu_j (dy)
&\leq&  \sum_{j=1}^N 4 \mu_j (K_j) \int_E u(x)^2 \mu_j (dx) \\
&\leq&     \sum_{j=1}^N 4 \mu_j (K_j) \| G_1 \mu_j\|_\infty \EE_1 (u, u).
\end{eqnarray*}
Thus there is a constant $C_0>0$ so that
$$
\EE_1 (u, u) \leq \EE^{(\lambda )} (u, u) \leq (1+C_0 \lambda) \EE_1 (u, u)
\quad \hbox{for every } u \in \FF.
$$
It follows that for every $\lambda > 0$, $(\EE^{(\lambda )}, \FF)$ is a regular Dirichlet form
on $L^2(E; m)$.

 \begin{thm}\label{T:4.2} For any increasing sequence $\{\lambda_n,n\geq 1\}$ of positive real numbers that increases to infinity, the Dirichlet form $(\mathcal{E}^{(\lambda_n)},\mathcal{F})$ is Mosco convergent to  the closed symmetric form $(\wt \EE,  \wt \FF)$ on $L^2(E; m)$.
\end{thm}

\noindent \begin{proof} Let $\{u_n,n\geq 1\}$ be a sequence in $L^2(E; m)$ that converges weakly to $u$ in $L^2(E; m)$ with
$\liminf_{n\rightarrow \infty}\mathcal{E}^{(\lambda_n)}(u_n,u_n)<\infty$. By taking a subsequence if necessary, we may and do assume that $\mathcal{E}^{(\lambda_n)}(u_n,u_n)$ converges, $\sup_{n\geq 1}\mathcal{E}_1^{(\lambda_n)}(u_n,u_n)<\infty$ and that the Cesaro mean sequence 
$\{v_n:=\sum_{k=1}^{n}u_k/n;n\geq 1\}$ is $\mathcal{E}^{(1)}_1$-convergent to some $v\in \mathcal{F}$. (The last property follows from Banach-Saks Theorem, see, for example, Theorem A.4.1 of \cite{CF}). As in particular, $v_n$
is $L^2(E; m)$-convergent to $v$, we must have $v=u$ $m$-a.e. on $E$. Hence $u\in \FF$ has a quasi-continuous version
which will still be denoted as $u$. Thus for every $k\geq 1$,
\begin{eqnarray}
\infty&>& \liminf_{n\rightarrow \infty}\mathcal{E}^{(\lambda_n)}(u_n,u_n)
\geq  \lim_{n\rightarrow \infty}
\mathcal{E}^{(\lambda_n)}(v_n,v_n)\geq \lim_{n\rightarrow \infty}\mathcal{E}^{(\lambda_k)}(v_n,v_n)=
\mathcal{E}^{(\lambda_k)}(u,u) \nonumber \\
\label{e:4.6}
&=& \mathcal{E}(u,u)+\lambda_k\sum_{j=1}^N\int_{K_j\times K_j}(u(x)-u(y))^2\mu_j(dx)\mu_j(dy).
\end{eqnarray}
Letting $k\rightarrow \infty$ in above inequality, we conclude that for each $j=1,\cdots,N,$
$$\int_{K_j\times K_j}(u(x)-u(y))^2\mu_j(dx)\mu_j(dy)=0.$$
This implies that $u$ is constant $\mu_j$-a.e. on $K_j$ and hence q.e. on $K_j$.
Thus  $u\in \wt \FF$ and by (\ref{e:4.6})
$$
\liminf_{n\rightarrow \infty}\mathcal{E}^{(\lambda_n)}(u_n,u_n)\geq \wt \EE (u,u),
$$
which establishes part (a) for the Mosco convergence.

To show part (b) of the Mosco convergence, it suffices to establish it for $u\in \wt \FF$ (for $u\neq
\wt \FF$, $\wt \EE (u,u)=\infty$ and so the property holds automatically).
Note that $\wt \FF\subset \mathcal{F}$. We take $u_n=u$ for every $n\geq 1$.  Then
$$
\mathcal{E}^{(\lambda_n)}(u_n,u_n) = \EE (u, u)= \wt \EE (u,u) \quad \hbox{for every } n\geq 1.
$$
This completes the proof of the theorem. \qed
\end{proof}

Let $X^n=\{X^n_t, t\geq 0; \P^n_x, x\in E\}$ be the Hunt process associated with the regular Dirichlet form
$(\EE^{(\lambda_n)}, \FF)$ on $L^2(E; m)$.
recall that $Y$ is the sticky Markov process with darning associated with the regular Dirichlet form $(\EE^*, \FF^*)$
on $L^2(E^*; \mu)$. 

\begin{thm}\label{T:4.3}
For every $0=t_0<t_1<\cdots t_k<\infty$ and bounded $\{f_j; 1\leq j\leq k\}\subset \overline{\wt \FF}$,
$$\lim_{n\to \infty} \E_m^n \left[ \prod_{j=0}^k f_j (X^n_{t_j}) \right] =
\E^*_\mu \left[ \prod_{j=0}^k f_j^* (Y_{t_j}) \right],
$$
where $f^*_j$ is defined by \eqref{e:4.4} with $f_j$ in place of $f$. 
\end{thm}

\pf For simplicity, we prove the theorem for $k=2$; the other cases are  similar.
Note that the semigroup $\{P^n_t; t\geq 0\}$ associated with the regular Dirichlet form
$(\EE^{(\lambda_n)}, \FF)$ on $L^2(E; m)$ is given by
$ P^n_t f(x)= \E^n_x [ f(X^n_t)]$, while, in view of Theorem \ref{T:4.1}, the semigroup
$\{\wh P_t; t\geq 0\}$ associated with the closed symmetric form $(\wt \EE, \wt \FF)$ on $\overline {\wt \FF}\subset L^2(E; m)$
is given by
$$
\wh P_t f = T^{-1} P^*_t f^*, \quad \hbox{where } P^*_t f^* (x)=\E_x \left[ f^* (Y_t)\right],
$$
for $f\in \overline {\wt \FF}$.
By Theorems \ref{T:4.2} and \ref{T:2.3}, $P^n_t f$ converges to $\wh P_t f$ for every $f\in \overline {\wt \FF}$
and $t>0$. It follows that $f_1P^n_{t_2-t_1}f_2 $ converges to $f_1\wh P_{t_2-t_1} f_2$ in $L^2(E; m)$.
Since $f_1\wh P_{t_2-t_1} f_2 \in \overline {\wt \FF}$, it follows
\begin{eqnarray*}
&& \lim_{n\to \infty} \|P^n_{t_1} ( f_1P^n_{t_2-t_1}f_2 ) -\wh P_{t_1} (f_1\wh P_{t_2-t_1} f_2) \|_{L^2(E; m)} \\
&\leq & \lim_{n\to \infty} \left( \|P^n_{t_1} ( f_1P^n_{t_2-t_1}f_2 -f_1\wh P_{t_2-t_1} f_2) \|_{L^2(E; m)}
+ \|  P^n_{t_1} (f_1\wh P_{t_2-t_1} f_2)- \wh P_{t_1} (f_1\wh P_{t_2-t_1} f_2) \|_{L^2(E; m)} \right) \\
&\leq & \lim_{n\to \infty}  \|  f_1P^n_{t_2-t_1}f_2 -f_1\wh P_{t_2-t_1} f_2 \|_{L^2(E; m)} =0.
\end{eqnarray*}
Hence we have
\begin{eqnarray*}
\lim_{n\to \infty} \E_m^n \left[ \prod_{j=0}^k f_j (X^n_{t_j}) \right]
&=& \lim_{n\to \infty} \int_E f_0(x) P^n_{t_1} ( f_1P^n_{t_2-t_1}f_2 ) (x) m(dx) \\
&=&  \int_E f_0(x) \wh P_{t_1} ( f_1 \wh P_{t_2-t_1}f_2 ) (x) m(dx) \\
&=&  \int_{E^*} f^*_0(x) P^*_{t_1} ( f^*_1 P^*_{t_2-t_1}f^* _2 ) (x) \mu (dx) \\
&=& \E^*_\mu \left[ \prod_{j=k}^n f_j^* (Y_{t_j}) \right] . 
\end{eqnarray*}
\qed

Theorem \ref{T:4.3} says that $X^n$ converges to the sticky Markov process with darning $Y$ in the
finite dimensional sense for all the testing functions that are constant on each $K_j$.

\medskip

When $(\EE, \FF)$ is a local Dirichlet form (or, equivalently, when $X$ is a diffusion on $E$)
 and each $K_j$ is connected and has positive measure, it is possible to approximate  sticky diffusions
 with darning by increasing the diffusion coefficients on each $K_j$ to infinity. 
This provides a very intuitive picture for shorting of each $K_j$ -- achieved by increasing
the conductance on $K_j$ to infinity.
 We illustrate this by the following example.

Suppose that $A(x)=(a_{ij}(x))_{1\leq i, j\leq d}$
is a matrix-valued function on $\R^d$ that is uniformly elliptic and bounded, 
and $\rho$ is a measurable function on 
$\R^d$ that is bounded between two positive constants.
Define  $\FF=W^{1, 2}(\R^d)=\left\{u\in L^2(\R^d; dx): \nabla u \in L^2 (\R^d; dx) \right\}$ and  
\begin{equation}\label{e:4.7}
\EE (u, v) = \frac12 \int_{\R^d} \sum_{i,j=1}^d a_{ij}(x) \frac{\partial u(x)}{\partial x_i} 
\frac{\partial v(x)}{\partial x_j} \rho (x) dx, \quad u, v \in W^{1,2}(\R^d).
\end{equation}
Then $(\EE, \FF)$ is a strongly local regular Dirichlet form on $L^2(\R^d; m)$, where $m(dx):= \rho (x) dx$.
It uniquely determines an $m$-symmetric diffusion process $X$ on $\R^d$ whose infinitesimal generator is 
$$
\LL=\frac1{2\rho (x)} \sum_{i,j=1}^d \frac{\partial}{\partial x_i} \left( \rho (x) a_{ij}(x) \frac{\partial}{\partial x_j}  \right).
$$
Let $\{K_j; 1\leq j\leq N\}$ be a finite number of disjoint compact sets which are the closure of
 non-empty connected open sets. Let $\{\lambda_n; n\geq 1\}$ 
 be an increasing sequence of positive numbers that increases to $\infty$. Define
$$ 
\EE^{(n)}(u, v) =\EE (u, v) + \lambda_n \sum_{l=1}^N \int_{K_l}
 \sum_{i,j=1}^d a_{ij}(x) \frac{\partial u(x)}{\partial x_i} 
\frac{\partial v(x)}{\partial x_j} \rho (x) dx, \quad u, v \in W^{1,2}(\R^d).
$$
Clearly for every $n\geq 1$,
$(\EE^{(n)}, \FF)$ is a regular $m$-symmetric strongly local Dirichlet form on $L^2(\R^n; m)$
and so there is an $m$-symmetric diffusion process $X^{(n)}$ associated with it. 
Let $\wt \FF$ be defined as in \eqref{e:4.1},  and $\wt \EE := \EE$ on $\wt \FF$.

Define $D=E\setminus \cup_{j=1}^N K_j$. We short (or collapse) each $K_j$ into a single point $a^*_j$.
Formally, by identifying each $K_j$ with a single point $a^*_j$, we
can get an induced topological space $E^*:=D\cup \{a_1^*, \dots,
a_N^*\}$ from $E$, with a neighborhood of each $a^*_j$ defined
as $(U\cap D)\cup \{a^*_j\}$ for some neighborhood $U$ of
$K_j$ in $E$. We define a measure $\mu$ on $E^*$ by setting $\mu=m$ and $D$ and  $m^* (a^*_j)=m(K_j)$.
Let $(\EE^*, \FF^*)$ be defined from $(\EE, \FF)$ as in \eqref{e:3.2}-\eqref{e:3.3}.
Then $(\EE^*, \FF^*)$ is a regular Dirichlet form on $L^2(E^*; \mu)$.
There is an associated diffusion process $Y$ on $E^*$, which we call sticky 
diffusion process with darning. If we take $m^*$ defined by $m^*(A):=\mu (A\cap D)$,
the diffusion process $X^*$ on $E^*$ associated with the regular Dirichlet form
$(\EE^*, \FF^*)$ on $L^2(E^*; m^*)$ is called diffusion process with darning.
These two processes are related to each other by a time change. 

\begin{thm}\label{T:4.4} Suppose $\{\lambda_n; n\geq 1\}$ is an increasing sequence of positive numbers that increases to infinity.  

{\rm (i)} The Dirichlet form $(\mathcal{E}^{(\lambda_n)},\mathcal{F})$ is Mosco convergent to  the closed symmetric form 
$(\wt \EE,  \wt \FF)$ on $L^2(E; m)$.

{\rm (ii)} Let $X^n=\{X^n_t, t\geq 0; \P^n_x, x\in E\}$ be the Hunt process associated with the regular Dirichlet form
$(\EE^{(\lambda_n)}, \FF)$ on $L^2(E; m)$. Then $X^n$ converges in the finite dimensional distribution
sense of Theorem \ref{T:4.3} to the sticky diffusion with darning $Y$ on $E^*$.
\end{thm}

\noindent \begin{proof} 
The proof is similar to that for Theorem \ref{T:4.2}.
For reader's convenience, we spell out the details.
 Let $\{u_n,n\geq 1\}$ be a sequence in $L^2(E; m)$ that converges weakly to $u$ in $L^2(E; m)$ with
$\liminf_{n\rightarrow \infty}\mathcal{E}^{(\lambda_n)}(u_n,u_n)<\infty$. By taking a subsequence if necessary, we may and do assume that $\mathcal{E}^{(\lambda_n)}(u_n,u_n)$ converges, $\sup_{n\geq 1}\mathcal{E}_1^{(\lambda_n)}(u_n,u_n)<\infty$ and that the Cesaro mean sequence 
$\{v_n:=\sum_{k=1}^{n}u_k/n;n\geq 1\}$ is $\mathcal{E}^{(1)}_1$-convergent to some $v\in \mathcal{F}$.   As in particular, $v_n$
is $L^2(E; m)$-convergent to $v$, we must have $v=u$ $m$-a.e. on $E$. Hence $u\in \FF$ has a quasi-continuous version
which will still be denoted as $u$. For every $k\geq 1$,
\begin{eqnarray}
\infty&>& \liminf_{n\rightarrow \infty}\mathcal{E}^{(\lambda_n)}(u_n,u_n)
\geq  \lim_{n\rightarrow \infty}
\mathcal{E}^{(\lambda_n)}(v_n,v_n)\geq \lim_{n\rightarrow \infty}\mathcal{E}^{(\lambda_k)}(v_n,v_n)=
\mathcal{E}^{(\lambda_k)}(u,u) \nonumber \\ 
&=& \mathcal{E}(u,u)+\lambda_k\sum_{j=1}^N \sum_{l=1}^N \int_{K_l}
\sum_{i,j=1}^d a_{ij}(x) \frac{\partial u(x)}{\partial x_i} 
\frac{\partial v(x)}{\partial x_j} \rho (x) dx .  \label{e:3.8}
\end{eqnarray}
Letting $k\rightarrow \infty$ in above inequality yields   
$\nabla u=0$ a.e. on $K_j$.
This implies that $u$ is constant a.e. in the interior of $K_j$.
Since $u$ is $\EE$-quasi-continuous on $\R^d$,  $u$ is constant $\EE$-q.e. on each $K_j$.
Hence $u\in \wt \FF$ and by \eqref{e:3.8}
$$
\liminf_{n\rightarrow \infty}\mathcal{E}^{(\lambda_n)}(u_n,u_n)\geq \wt \EE (u,u),
$$
which establishes part (a) for the Mosco convergence.

To show part (b) of the Mosco convergence, it suffices to establish it for $u\in \wt \FF$.
Note that $\wt \FF\subset \mathcal{F}$. We take $u_n=u$.  Then
$$
\mathcal{E}^{(\lambda_n)}(u_n,u_n) = \EE (u, u)= \wt \EE (u,u) \quad \hbox{for every } n\geq 1.
$$
This completes the proof that the Dirichlet form $(\mathcal{E}^{(\lambda_n)},\mathcal{F})$ is Mosco convergent to 
the closed symmetric form $(\wt \EE,  \wt \FF)$ on $L^2(E; m)$.

(ii) The proof is exactly the same as that for Theorem \ref{T:4.3}. 
\qed
\end{proof}

\section{Brownian motion on spaces with varying dimension}\label{S:5} 

A simple example of spaces with varying dimension is a large square with a thin
flag pole. Mathematically, it can be modeled by a plane with a vertical line installed 
on it:
\begin{equation}\label{e:5.1}
\R^2 \cup \R_+: =\left\{(x_1, x_2, x_3)\in \R^3: x_3=0 \textrm{ or } x_1=x_2=0 \hbox{ and } x_3>0\right\}. 
\end{equation}
Spaces with varying dimension arise  in  many disciplines including statistics, physics  and engineering 
(e.g.  molecular dynamics, plasma dynamics).
It is natural to study Brownian motion and ``Laplacian operator" on such spaces. 
Intuitively, Brownian motion on space $\R^2\cup \R_+$ should behave like a two-dimensional
Brownian motion when it is on the plane, and like a one-dimensional Brownian motion
when it is on the vertical line (flag pole).  However  the space $\R^2\cup \R$  is quite singular in the sense
that the base $O$ of the flag pole where the plane and the vertical line meet 
is a singleton. A singleton
would never be visited by a two-dimensional Brownian motion, which means 
 Brownian motion starting from a point on the plane will never visit $O$.
Hence there is no chance for such a process to climb up the flag pole. 
The solution is to 
 collapse or short (imagine putting
an infinite conductance on) a small closed disk $\overline{ B(0, \eps)}\subset \R^2$ centered at the origin into a point $a^*$
 and   consider the resulting
Brownian motion with darning on the collapsed plane,
for which $a^*$ will be visited. Through $a^*$ a vertical  pole can be installed 
and  one can construct Brownian motion
with varying dimension (BMVD) on $\R^2 \cup \R_+$ by joining together
the Brownian motion with
darning on the plane and the one-dimensional Brownian motion along the pole.
This is done in \cite{CL} through a Dirichlet form method. 

To be more precise, the state space of BMVD  is defined as follows.
Fix $\eps>0$ and $p>0$.
 Let $D_0=\R^2\setminus  {\overline B (0, \eps)} $. By identifying the closed ball
$\overline {B(0, \eps)}$ with a singleton denoted by $a^*$, we can introduce a topological space $E^*:=D_0\cup \{a^*\}\cup \R_+$, with the origin of $\R_+$ identified with $a^*$ and with the topology on
$E^*$ induced from that of $\R^2 \cup \R_+$. Let $m^*_p$ be the measure on $E^*$ whose restriction on $\R_+$ and $D_0$ is the Lebesgue measure multiplied by $p$ and $1$, respectively. 

\medskip

\begin{defn}\label{D:5.1} \rm  Let $\eps >0$ and $p>0$. A Brownian motion with varying dimension
 (BMVD in abbreviation) on $E^*$ with parameters $(\eps, p)$ on $E^*$  is an $m^*_p$-symmetric diffusion $X^*$ on $E$ such that
\begin{description}
\item{(i)} its part process in $\R_+$ or $D_0$ has the same law as standard Brownian motion killed upon leaving $\R_+$ or $D_0$,
respectively;
\item{(ii)} it admits no killings on $a^*$.
\end{description}
\end{defn}

It follows from the $m^*_p$-symmetry of $X^*$ and the fact  $m^*_p(\{a^*\})=0$
that BMVD $X$ spends zero Lebesgue amount of time at $a^*$.
It is shown in \cite[Theorem 1.2]{CL} that 
 for every $\eps>0$ and $p>0$, BMVD with parameters $(\eps, p)$
exists and is unique in law.  
In fact,  BMVD on $E^*$ can be constructed as the $m^*_p$-symmetric Hunt process
 associated with the regular Dirichlet form $(\EE^*, \FF^*)$ on $L^2(E^*; m^*_p)$  given by
\begin{eqnarray}
\FF^* &= &  \left\{f: f|_{D_0}\in W^{1,2}(D_0),  \, f|_{\R_+}\in W^{1,2}(\R_+),
\hbox{ and }
f (x) =f (0) \hbox{ q.e. on } \partial D_0 \right\},   \label{e:5.2}
\\
\EE^* (f,g) &=& \frac12 \int_{D_0}\nabla f(x) \cdot \nabla g(x) dx+\frac{p}2\int_{\R_+}f'(x)g'(x)dx .  \label{e:5.3}
\end{eqnarray}
Here for an open set $U\subset \R^d$, $W^{1,2}(U)$ is the Sobolev space on $U$ of order $(1, 2)$; that is,
$$ W^{1,2}(U)= \left\{  f\in L^2(U; dx): \nabla f \in L^2(U; dx)  \right\}.
$$
Sample path properties  of $X^*$ including that at the base point and the two-sided transition density function estimates
have been studied in \cite{CL}.  Roughly speaking, when BMVD $X^*$ is at the base point $a^*$,  it enters the pole with probability
$\frac{p}{2\pi \eps +p}$ and enters the punched plane $D_0$ with probability $\frac{2\pi \eps}{2\pi \eps + p}$; see \cite[Proposition 4.3]{CL}. 

\medskip

We will show in this section that BMVD on $E^*$ can be approximated
by Brownian motion in the plane with a vertical cylinder whose horizontal motion
on the cylinder is a circular Brownian motion moving at fast speed.  
Let 
\begin{equation}\label{e:5.4}
 E: =\left\{(x_1, x_2, x_3)\in \R^3:  x_1^2+x_2^2\geq \eps^2 \hbox { and }  x_3=0 \hbox{ or } 
x_1^2+x_2^2=\eps^2  \hbox{ and } x_3>0\right\} .
\end{equation}
That is, $E$ is $D_0\times \{0\}$ with a vertical cylinder with base radius $\eps$ sitting on top of $(\partial D_0) \times \{0\}$. 
Let $\overline m_p$ be the measure on $E$ whose restriction on $D_0 \times \{0\}$ is the two-dimensional Lebesgue measure
and its restriction to the cylinder $\partial D_0 \times [0, \infty)$ is the  Lebesgue surface measure multiplied by $p/(2\pi \eps)$. 
When there is no danger of confusion,   we identify  $D_0$ with  $D_0\times \{0\}$. 
The space $E$ is a two-dimensional Lipschitz manifold. For every $\lambda >0$,
we can run an $\overline m_p$-symmetric diffusion $X^{(\lambda)}$ on $E$ that behaves as Brownian motion on $D_0\times \{0\}$
and behaves like  $(B_{\lambda t}, W_t)$   while on the cylinder $\partial D_0 \times (0, \infty)$.
Here $B_t$ is a  standard circular Brownian motion on  $\partial D_0$ and $W_t$ is a one-dimensional Brownian motion.    
We will show that as $\lambda \to \infty$, $X^{(\lambda)}$ converges in the finite-dimensional distribution sense,
after a suitable identification,  to the BMVD $X^* $ on $E^*$; see Theorem \ref{T:5.3} for a precise statement. Note that the state space $E$ of $X^{(\lambda)}$,
which is of dimension two, is different from the state space $E^*$ of $X^*$ with varying dimension. The space $E^*$ can be viewed as $E$  with  
the cylinder $\partial D_0 \times [0, \infty)$ collapsed into one single half line $\{a^* \}\times [0, \infty)$. 
The main difference,  when compared with Brownian motion with darning in Section \ref{S:3},  is that   here we collapse every circle
$\partial D_0 \times \{ z\}$ into one point $\{a^*\}\times \{ z\}$ and there are a continuum of such circles to collapse.  
However the ideas developed in Section \ref{S:4} can be adapted to establish the convergence of $X^{(\lambda)}$ to  BMVD $X^*$,
and we spell out the details in what follows.  

First we give a precise construction of $X^{(\lambda)}$ via a Dirichlet form approach.
First, we introduce Sobolve space $W^{1,2}(E)$ of order $(1, 2)$ on $E$. 
For convenience, let $S:=\partial D_0 \times (0, \infty)$. 
The space $S$ can be identified with the infinite rectangle $[0, 2\pi \eps]\times (0, \infty)$ with 
points $(0, z)$ and $(2\pi \eps, z)$ identified. We denote the pull back measure on $S$ of the Lebesgue measure on $ [0, 2\pi \eps]\times (0, \infty)$ by $m$. 
Functions on $S$ can be parametrized 
by $(t, z)\in [0, 2\pi \eps)\times [0, \infty)$. We define
$$
W^{1,2}(S) =\left\{f=f (x, z)\in L^2 (S; m): \left| \partial_t f\right|
+\left| \partial_z f  \right| \in L^2(S; m) \right\}.
$$
Note that $W^{1,2}(D_0)$ and $W^{1,2}(S)$ are  the Dirichlet spaces for the reflecting Brownian motion
on $D_0$ and on the cylinder $S$, respectively. So for every $f\in W^{1, 2}(D_0)$
and $g\in W^{1,2}(S)$,  their quasi-continuous versions are well defined on $\partial D_0=\partial S$
quasi-everywhere, which we call the trace on the circle $\partial D_0$ and we denote them by
$f|_{\partial D_0}$ and $g|_{\partial D_0}$, respectively.
Define
\begin{eqnarray*}
 W^{1,2}(E):=\left\{ f: f_1:=f |_{D_0} \in W^{1,2}(D_0), f_2=f|_S \in W^{1,2}(S) \hbox{ and }
 f_1 |_{\partial D_0} = f_2|_{\partial D_0}   \hbox{ q.e.} \right\}
\end{eqnarray*}
For $f\in W^{1,2}(E)$, define its norm $\|f \|_{1,2}$ by
$$
\| f\|_{1,2}^2 = \int_{D_0}  | \nabla f(x)|^2 dx + \int_{S} \left( |\partial_t f(t, z)|^2 + 
| \partial_z f (t, z)^2 \right) m(dt dz).
$$
It is easy to see that $W^{1,2} (E)$ is the $\| \cdot \|_{1,2}$-completion of the following subspace
of continuous functions on $E$:
$$
\left\{ f\in C(E): \  f |_{D_0} \in W^{1,2}(D_0),  \ f|_S \in W^{1,2}(S)    \right\}
$$
Now for every $\lambda >0$, define $\FF^{(\lambda)}=W^{1,2}(E)$ and for $f\in \FF^{(\lambda)}$, 
\begin{eqnarray} \label{e:5.5} 
\EE^{(\lambda)} (f, f ) &=& \frac12 \int_{D_0} |\nabla f(x)|^2 dx 
+ \frac{\lambda p}{4\pi \eps} \int_0^\infty \left(\int_0^{2\pi \eps} |\partial_t f(t, z)|^2  dt\right) dz  \nonumber \\
&& + \frac{p}{4\pi \eps}   \int_0^{2\pi \eps} \left( \int_0^\infty |\partial_z f(t, z)|^2  dz \right) dt.
\end{eqnarray}
The last two terms in the right hand side of \eqref{e:5.5} represents the $\EE^{(\lambda)}$-energy of $f$ on the cylinder $S$. 
It is easy to check that $(\EE^{(\lambda)}, \FF^{(\lambda)})$ is a symmetric regular strongly local Dirichlet form on $L^2(E; \overline m_p)$ and so it uniquely determines a symmetric Hunt process $X^{(\lambda)}$ on $E$. 
Using the part Dirichlet form of $(\EE^{(\lambda)}, \FF^{(\lambda)})$ on $D_0$ and $S$, respectively, 
it is easy to see \cite{CF, FOT} that the part process of $X^{(\lambda)}$ in $D_0$ and $S$ are 
the part process of two dimension Brownian motion in $D_0$ and $(  B_{\lambda t}, W_t)$ on $S$, respectively.
Here $B_t$ is the circular Brownian motion on $\partial D_0$ and $W_t$ is Brownian motion on $(0, \infty)$ independent 
of $B_t$.

Let
\begin{equation}\label{e:5.6}
\wt \FF= \{ f\in W^{1,2}(E): f =\hbox{constant $\EE$-q.e. on } \partial D_0\times \{z\} \hbox{ for each } z\geq 0 \} .
\end{equation}
Note that since each circle $\partial D_0\times \{z\} $ is of positive $\EE^{(\lambda)}$-capacity for every $\lambda >0$, 
$(\EE^{(\lambda)}, \wt \FF)$ is a closed symmetric  Markovian bilinear form on $L^2(E; \overline m_p)$ but $\wt \FF$ is not
dense in $L^2(E; m)$. Note that $\EE^{(\lambda)}=\EE^{(1)}$ on $\wt \FF$ for every $\lambda >0$.
To emphasize its dependence on the domain of definition, we write $(\wt \EE, \wt \FF)$ for $(\EE^{(1)}, \wt \FF)$.
Denote by $\Pi$ the orthogonal projection of $L^2(E; \overline m_p)$ onto the closure $\overline {\wt \FF}$ 
of $ \wt \FF$ in $L^2(E; \overline m_p)$.
Let $\{\wt G_\alpha, \alpha>0\}$ be the resolvent associated with
$(\wt \EE, \wt \FF)$ on $L^2(E; m)$, and $\{\wh P_t; t\geq 0\}$ and $(\wh G_\alpha; \alpha >0\}$
the strongly continuous semigroup and resolvent of the closed symmetric form  $(\wt \EE, \wt \FF)$ on $\overline {\wt \FF}$, respectively.
We know from \eqref{e:2.3a} that $\wt G_\alpha f = \wh G_\alpha (\Pi f)$.
 We next identify $\wh P_t$ and $\wh G_\alpha$, as well as $\Pi$.

The following map $T$ establishes a one-to-one and onto correspondence between
$(\wt \EE, \wt \FF)$ on $\overline {\wt \FF}\subset L^2(E; \overline m_P)$
and $(\EE^*, \FF^*)$ on $L^2(E^*; m^*_p)$: for every $f\in \wt \FF$,
\begin{equation}\label{e:5.7} 
Tf=f \ \hbox{ on } D_0 \quad \hbox{ and } \quad Tf ( a^*, z )=f( \partial D_0 \times \{z\}) \ \hbox{ for } z\geq 0.
\end{equation}
(For $g\in \FF^*$, $T^{-1} g(x)=g(x)$ for $x\in D_0$ and $T^{-1} g(t, z)= g(a^*, z)$  for $(t, z)\in S$.)
The map $T$ has the property that for every $f\in \wt \FF$,
\begin{equation}\label{e:5.8}
 \wt \EE (f, f)= \EE^*(Tf, Tf) \quad \hbox{and} \quad \| f\|_{L^2(E;  \overline m_p)}=\| Tf \|_{L^2(E^*; m^*_p)}.
\end{equation}
In other words, $T$ is an isometry between $(\wt \EE, \wt \FF)$ on $\overline {\wt \FF}\subset L^2(E; \overline m_p)$
 and $(\EE^*, \FF^*)$ on $L^2(E^*;   m^*_p)$
both in $\EE$ and in the $L^2$ sense.
Denote by $\{ G^*_\alpha; \alpha >0\}$ and $\{P^*_t; t\geq 0\}$ the resolvent and semigroup
associated with the regular Dirichlet form
$( \EE^*, \FF^*)$ on $L^2(E^*; m^*_p)$, and $\{ \wh G_\alpha; \alpha >0\}$ and $\{\wh P_t; t\geq 0\}$
 the resolvent and semigroup on $\overline {\wt \FF}$
associated with the closed symmetric  form
$(\wt \EE, \wt \FF)$ on $L^2(E; \overline m_p)$.

For every $f\in L^2(E; \overline m_p)$, we can define a function $f^*$ on $E^*$ by setting 
$f^*=f$ on $D_0$ and
\begin{equation}\label{e:5.9}
f^* (a^*, z)=   \frac1{2 \pi \eps} \int_{\partial D_0} f(t, z) dt   
\quad \hbox{for } z>0.
\end{equation} 
Note that by Fubini theorem, $f^*(a^*, z)$ is well defined for a.e. $z\in (0, \infty)$. 
For an $\EE^*$-quasi-continuous function $g$ defined on $E^*$, we define
$$
T^{-1} g = g(x)  \quad \hbox{on } D_0 \qquad \hbox{and} \qquad 
T^{-1} (t, z) = g( a^*, z ) \quad \hbox{for }  (t, z)\in \partial D_0 \times [0, \infty) . 
$$
Clearly,  
\begin{equation}\label{e:5.10}
( T^{-1} g )^*= g \qquad \hbox{and} \qquad \| g \|_{L^2(E^*; m^*_p)}=\| T^{-1} g\|_{L^2(E; \overline m_p)}  .
\end{equation}
Since $T$ extends to be an isometry between $\overline{\wt \FF}\subset L^2(E; \overline m_p)$ and $L^2(E^*; m^*_p)$, 
the above defined $T^{-1}$ extends to be an isometry between $L^2(E^*; m^*_p)$ and $\overline{\wt \FF}\subset L^2(E; \overline m_p)$.

\begin{thm}\label{T:5.2} \begin{description}
\item{\rm (i)}
For $f\in L^2(E; \overline m_p)$, $\Pi f = T^{-1} f^*$ $\overline m_p$-a.e. on $E$.

\item{\rm (ii)}  $\wh G_\alpha f =T^{-1} G^*_\alpha f^*$ for $f\in L^2(E; \overline m_p)$.
\end{description}

\item{\rm (iii)} For $f\in \overline {\wt \FF}$,  $t>0$ and $\alpha>0$,
\begin{equation}\label{e:5.11}
\wh P_t f=  T^{-1} P^*_t (f^*) \quad \hbox{and} \quad  \wh G_\alpha f =
 T^{-1} G^*_\alpha (f^*).
\end{equation}
\end{thm}

\pf (i) Let ${\mathcal C}=\FF^*\cap C_c (E^*)$, which is a core of the regular Dirichlet form
$(\EE^*, \FF^*)$ on $L^2(E^* ; m^*_p)$. In particuar,  ${\mathcal C}$ is $L^2$-dense in $L^2(E^*;  m^*_p)$.
It  follows from \eqref{e:5.8} and \eqref{e:5.10} that $T^{-1}{\mathcal C}$ is $L^2(E; \overline m_p)$-dense in $\overline {\wt \FF}$
and by Fubini's theorem, 
$$
\overline {\wt \FF}=  \left\{f\in L^2(E; \overline m_p): f \hbox{ is constant a.e. on } \partial D_0 \times \{z\} \hbox{ for a.e. } z>0 \right\}.
$$
Thus for every $f\in L^2(E; \overline m_p)$, $ f^*\in \overline {\wt \FF}$.
On the other hand, it is clear that for $f\in L^2(E; \overline m_p)$,
$$
(f, T^{-1}g)_{L^2(E; \overline m_p)} = (T^{-1} f^*, T^{-1}g)_{L^2(E; \overline m_p)}
\quad \hbox{for every } g\in {\mathcal C}.
$$
Thus we conclude that  $\Pi f =T^{-1} f^*$ $\overline m_p$-a.e.

(ii) Let $f\in L^2(E;  \overline m_p)$.
For every $\alpha >0$ and $g\in \wt \FF$, it follows from \eqref{e:5.8} that
\begin{eqnarray*}
\wt \EE_\alpha ( \wt G_\alpha f, g)
&=& \int_E f(x) g(x)  \overline m_P (dx) = \int_{E^*} f^*(x) T g(x) m^*_p (dx) \\
&=& \EE^* ( G^*_\alpha f^*, Tg)+\alpha \int_{E^*} G^* f^* (x) Tg (x) m^*_p  (dx) \\
&=& \wt  \EE_\alpha (T^{-1} G^*_\alpha f^*, g).
\end{eqnarray*}
We thus have $\wt G_\alpha f =T^{-1} G^*_\alpha f^*$.

(iii) This follows immediately from (i), (ii) and \eqref{e:2.3a} that for $f\in \overline {\wt \FF}$,
$$
\wh G_\alpha f =\wt G_\alpha f = T^{-1} G^*_\alpha f^*.
$$
It is clear that $T_t f:= T^{-1}P^*_t f^*$ defines a symmetric strongly continuous contraction semigroup
on $\overline {\wt \FF}\subset L^2(E; \overline m_p)$, as
 $\{P^*_t; t\geq 0\}$ is a strongly continuous contraction semigroup on $L^2(E^*; m^*_p)$.
 Moreover, for every $\alpha>0$, $\int_0^\infty e^{-\alpha t} T_t f dt = T^{-1} G^*_\alpha f^*=\wh G_\alpha f$.
 We thus conclude that $T_t =\wh P_t$.
\qed

\begin{thm}\label{T:5.3} Suppose $\{\lambda_n; n\geq 1\}$ is an increasing sequence of positive numbers that increases to infinity.  

{\rm (i)} The Dirichlet form $(\mathcal{E}^{(\lambda_n)},\mathcal{F})$ is Mosco convergent to the closed symmetric form 
$(\wt \EE,  \wt \FF)$ on $L^2(E; \overline m_p)$.

{\rm (ii)} Let $X^n=\{X^n_t, t\geq 0; \P^n_x, x\in E\}$ be the Hunt process associated with the regular Dirichlet form
$(\EE^{(\lambda_n)}, \FF)$ on $L^2(E;  \overline m_p)$. Then $X^n$ converges in the finite dimensional distribution
sense of Theorem \ref{T:4.3} to the BMVD $X^*$ on $E^*$. 
\end{thm}

\noindent \begin{proof}  
  Let $\{u_n,n\geq 1\}$ be a sequence in $L^2(E; \overline m_p)$ that converges weakly to $u$ in $L^2(E; \overline m_p)$ with
$\liminf_{n\rightarrow \infty}\mathcal{E}^{(\lambda_n)}(u_n,u_n)<\infty$. By taking a subsequence if necessary, we may and do assume that $\mathcal{E}^{(\lambda_n)}(u_n,u_n)$ converges, $\sup_{n\geq 1}\mathcal{E}_1^{(\lambda_n)}(u_n,u_n)<\infty$ and that the Cesaro mean sequence 
$\{v_n:=\sum_{k=1}^{n}u_k/n;n\geq 1\}$ is $\mathcal{E}^{(1)}_1$-convergent to some $v\in \mathcal{F}$.   Since  $v_n$
is $L^2(E; \overline m_p)$-convergent to $v$, we must have $v=u\,$ $\overline m_p$-a.e. on $E$. Hence $u$ has a quasi-continuous version
which will still be denoted as $u$. Thus for every $k\geq 1$,
\begin{eqnarray}
\infty&>& \lim_{n\rightarrow \infty}\mathcal{E}^{(\lambda_n)}(u_n,u_n)
\geq  \lim_{n\rightarrow \infty}
\mathcal{E}^{(\lambda_n)}(v_n,v_n)\geq \lim_{n\rightarrow \infty}\mathcal{E}^{(\lambda_k)}(v_n,v_n)=
\mathcal{E}^{(\lambda_k)}(u,u) \nonumber \\ 
&=&  \frac12 \int_{D_0} |\nabla u(x)|^2 dx 
+ \frac{\lambda_k p}{4\pi \eps} \int_0^\infty \left(\int_0^{2\pi \eps} |\partial_t u(t, z)|^2  dt\right) dz  \nonumber \\
&&  + \frac{p}{4\pi \eps}   \int_0^{2\pi \eps} \left( \int_0^\infty |\partial_z u(t, z)|^2  dz \right) dt.  \label{e:5.12}
\end{eqnarray}
Letting $k\rightarrow \infty$ in above inequality, we conclude that  there is a subset ${\cal N}\subset (0, \infty)$
having zero Lebesgue measure 
so that for every $z\in (0, \infty) \setminus {\cal N}$, 
$ \int_0^{2\pi \eps} |\partial_t u(t, z)|^2  dt =0 $. 
This implies that for every $z \in (0, \infty) \setminus {\cal N}$, $t\mapsto u(t,z)$ is equals to a constant $u(z)$   a.e.  and hence   $\EE^{(1)}$-q.e.  on $ [0, 2\pi)$.
For $0<z_1<z_2$ in $(0, \infty) \setminus {\cal N}$, by Cauchy-Schwartz inequality, 
\begin{eqnarray*} 
|u(z_2)-u(z_1)| &=& \frac1{2\pi \eps} \left| \int_0^{2\pi \eps} (u(t, z_2)-u(t, z_1))dt \right|
= \frac1{2\pi \eps} \left| \int_0^{2\pi \eps} \int_{z_1}^{z_2} \partial_z u(t, z) dz dt \right|  \\
&\leq&  \frac1{\sqrt{2\pi \eps} } \left( \int_0^{2\pi \eps} \int_{z_1}^{z_2} |\partial_z u(t, z)|^2 dz dt \right) |z_2-z_1|^{1/2}.  
\end{eqnarray*}
This shows that $u(z)$ is a H\"older continuous function on $[0, \infty)$.
Since each horizontal circle and each vertical line on the cylinder is of positive $\EE^{(1)}$ capacity
and $u$ is $\EE^{(1)}$-quasi-continuous on $E$, it follows that an $\EE^{(1)}$-quasi-continuous version 
of $u$ can be taken so that $u(t, z)= u(z)$ for every $z\geq 0$ and $t\in [0, 2\pi)$
(such defined function is continuous on the cylinder $S$). 
Hence $u\in \wt \FF$ and by \eqref{e:5.12}
$$
\liminf_{n\rightarrow \infty}\mathcal{E}^{(\lambda_n)}(u_n,u_n)\geq \wt \EE (u,u),
$$
which establishes (a) for the Mosco convergence.

To show (b) of the Mosco convergence, it suffices to establish it for $u\in \wt \FF$.
Note that $\wt \FF\subset \mathcal{F}^{(\lambda)}$ for every $\lambda >0$. We take $u_n=u$.  Then
$$
\mathcal{E}^{(\lambda_n)}(u_n,u_n) = \EE (u, u)= \wt \EE (u,u) \quad \hbox{for every } n\geq 1.
$$
This proves that the Dirichlet form $(\mathcal{E}^{(\lambda_n)},\mathcal{F})$ is Mosco convergent to $(\wt \EE,  \wt \FF)$ on $L^2(E; \overline m_p)$.

(ii) The proof is similar to    that for Theorem \ref{T:4.3} except using Theorem \ref{T:5.2} instead of  Theorem \ref{T:4.1}.
We omit its details here.  
\qed
\end{proof}

We remark that, since each horizontal circle on the cylinder that is to be collapsed into one single point has zero $\overline m_p$ measure,
so the limiting process of $X^n$ is just the BMVD $X^*$ on $E^*$, not a sticky one. 

For other related work and approaches on
 Markov processes living on spaces with
possibly different dimensions, we refer the reader to \cite{ES,  Ha, Ku}  and the references therein.

\section{Examples}\label{S:6}

In this section, we give some examples of the Dirichlet forms $(\EE, \FF)$, or equivalently symmetric Markov processes, 
for which the main results in Section \ref{S:4} are applicable.

\begin{example} \rm (Sticky diffusion process with darning)
Let $(\EE, \FF)$ be the strong local regular Dirichlet form on $L^2(\R^d; m)$ defined by \eqref{e:4.7},
where $m(dx)=\rho (x) dx$. Suppose that  $K_1, \dots, K_N$ are separated, non-$\EE$-polar compact 
(possibly disconnected) subsets of $E$.
Let $F=\cup_{j=1}^N K_j$ and $D=\R^d\setminus F$. We short (or collapse) each $K_j$ into a single point $a^*_j$.
By identifying each $K_j$ with a single point $a^*_j$, we
can get an induced topological space $E^*:=D\cup \{a_1^*, \dots,
a_N^*\}$ from $E$, with a neighborhood of each $a^*_j$ defined
as $(U\cap D)\cup \{a^*_j\}$ for some neighborhood $U$ of
$K_j$ in $E$. Let $\mu=m$ and $D$ and $\mu (a^*_j)= m(K_j)$.
Let $(\EE^*, \FF^*)$ be defined as in \eqref{e:3.2}-\eqref{e:3.3}.
Then it is a regular Dirichlet form on $L^2(E^*; \mu)$.
There is a unique   
diffusion process $Y$ on $E^*$ associated with it, which we call
sticky diffusion process with darning. When $(a_{ij}(x))\equiv I$, the identity
matrix, and $\rho \equiv 1$, $Y$ is called sticky Brownian motion with darning. 
For each $1\leq j\leq N$, take a finite smooth $\mu_j$ whose quasi-support is $K_j$ and having bounded 1-potential $G_1 \mu_j$. For each $\lambda >0$, let $\EE^{(\lambda)}$ be defined by \eqref{e:1.3}.
$(\EE^{(\lambda)}, \FF)$ is a regular Dirichlet form on $L^2(\R^d; m)$ and it determines a diffusion
process with jumps $X^{(\lambda)}$. By Theorem \ref{T:4.3}, for any increasing sequence $\{\lambda_n; n\geq 1\}$ that increases to infinity, $X^{(\lambda_n)}$ converges in the finite dimensional distribution in the sense of Theorem \ref{T:4.3} to the sticky diffusion process with darning $Y$ 
on $E^*$. 
\end{example}

\begin{example} \rm  (Sticky stable process with darning)
Suppose the metric measure space $(E, \rho, m)$ is a $d$-set; that is, there are positive
constants $c_1, c_2$ so that 
$$ 
c_1 r^d \leq m(B(x, r)) \leq c_2 r^d \quad \hbox{for every } x \in E
\hbox{ and } 0<r<1.
$$
Here $B(x, r):=\{y\in E: \rho (y, x)<r\}$ is the open ball centered at $x$ with radius $r$.
Suppose $c(x, y)$ is a symmetric function on $E\times E$ that is bounded between two positive constants, and $0<\alpha <2$. Define 
$$ 
\EE (f, f)=\int_{E\times E} (f(x)-f(y))^2 \frac{c(x, y)}{\rho (x, y)^{d+\alpha}} m(dx) m(dy) ,
$$
and  let $\FF$ be the closure of Lipschitz functions on $E$ with compact support under 
$\EE_1$, where $\EE_1 (f, f):= \EE (f, f)+ \int_E f(x)^2 m(dx)$. 
The bilinear form $(\EE, \FF)$ is a regular Dirichlet form on $L^2(E; m)$.
Its associated Hunt process $X$ is called $\alpha$-stable-like process on $E$ (cf. \cite{CK, CKW}).
 Suppose that  $K_1, \dots, K_N$ are separated, non-$\EE$-polar compact 
(possibly disconnected) subsets of $E$.
Let $F=\cup_{j=1}^N K_j$ and $D=E\setminus F$. We short (or collapse) each $K_j$ into a single point $a^*_j$.By identifying each $K_j$ with a single point $a^*_j$, we
can get an induced topological space $E^*:=D\cup \{a_1^*, \dots,
a_N^*\}$ from $E$, with a neighborhood of each $a^*_j$ defined
as $(U\cap D)\cup \{a^*_j\}$ for some neighborhood $U$ of
$K_j$ in $E$. Let $\mu=m$ and $D$ and $\mu (a^*_j)= m(K_j)$.
Let $(\EE^*, \FF^*)$ be defined as in \eqref{e:3.2}-\eqref{e:3.3}.
Then it is a regular Dirichlet form on $L^2(E^*; \mu)$.
There is a unique   
Hunt process $Y$ on $E^*$ associated with it, which we call
sticky $\alpha$-stable-like process with darning.  
For each $\lambda >0$, let $\EE^{(\lambda)}$ be defined by \eqref{e:1.3}.
$(\EE^{(\lambda)}, \FF)$ is a regular Dirichlet form on $L^2(\R^d; m)$ and it determines a jump diffusion $X^{(\lambda)}$. By Theorem \ref{T:4.3}, for any increasing sequence $\{\lambda_n; n\geq 1\}$ that increases to infinity, $X^{(\lambda_n)}$ converges in the finite dimensional distribution in the sense of Theorem \ref{T:4.3} to the sticky $\alpha$-stable-like process with darning $Y$ on $E^*$. 
\end{example}

Similarly, we can consider darning of symmetric diffusions with jumps studied in \cite{CKK} and their approximation by introducing 
additional jumps over the hulls $K_j$.

\vskip 0.3truein

\noindent {\bf Zhen-Qing Chen}

\smallskip \noindent
Department of Mathematics, University of Washington, Seattle,
WA 98195, USA

\noindent
Email: \texttt{zqchen@uw.edu}

\bigskip

\noindent {\bf Jun Peng}

\smallskip \noindent
School of Mathematics and Statistics, 
Central South University,  
Changsha, Hunan, 410075, China

\noindent
Email: \texttt{pengjun0825@163.com}

\medskip

\end{document}